 \documentclass[12pt]{article}

 \usepackage{mathrsfs,amsfonts}
 \usepackage{color}
 \usepackage{extarrows}
 \usepackage{amsmath}
\usepackage{color}
\usepackage[round]{natbib}
 \newtheorem{theorem}{Theorem}[section]
 \newtheorem{lemma}[theorem]{Lemma}
 \newtheorem{corol}[theorem]{Corollary}
 \newtheorem{prop}[theorem]{Proposition}
 
 \newtheorem{con1}[theorem]{Condition}
 \newtheorem{remark}[theorem]{Remark}

 \newtheorem{example}{Example}[section]

 \def\blemma{\begin{lemma}\sl{}\def\elemma{\end{lemma}}}
 \def\bproposition{\begin{prop}\sl{}\def\eproposition{\end{prop}}}
 \def\btheorem{\begin{theorem}\sl{}\def\etheorem{\end{theorem}}}
 \def\bcorollary{\begin{corol}\sl{}\def\ecorollary{\end{corol}}}
 \def\bexample{\begin{example}\rm{}\def\eexample{\end{example}}}

 \def\beqlb{\begin{eqnarray}}\def\eeqlb{\end{eqnarray}}
 \def\beqnn{\begin{eqnarray*}}\def\eeqnn{\end{eqnarray*}}
 
 \def\qqquad{\qquad\qquad}

 \def\benumerate{\begin{enumerate}}\def\eenumerate{\end{enumerate}}
 \def\bitemize{\begin{itemize}}\def\eitemize{\end{itemize}}\def\itm{\item}

 \def\proof{\noindent{\it Proof.~~}}\def\qed{\hfill$\Box$\medskip}

 \def\<{\langle}\def\>{\rangle}

 \def\mcr{\mathscr}\def\mbb{\mathbb}\def\mbf{\mathbf}

 \def\ar{\!\!&}

\begin{document}

\bigskip\bigskip

\centerline{\huge\bf A continuous-state polynomial}

\smallskip

\centerline{\huge\bf branching process}

\bigskip

\centerline{\Large Pei-Sen Li}

\bigskip

 \centerline{Institute for Mathematical Sciences, Renmin University of China}

\centerline{Beijing, 100872, China}

\centerline{Department of Mathematics and Statistics, Concordia University}

\centerline{Montreal, H3G 1M8, Canada}

\centerline{E-mail: \tt peisenli@mail.bnu.edu.cn}

\bigskip


\bigskip

\noindent\textbf{Abstract} 

A continuous-state polynomial branching process is constructed as the pathwise unique solution of a stochastic integral equation with absorbing boundary condition. The process can also be obtained from a spectrally positive L\'{e}vy process through Lamperti type transformations. The extinction and explosion probabilities and the mean extinction and explosion times are computed explicitly.  Some of those are also new for the classical linear branching process. We present necessary and sufficient conditions for the process to extinguish or explode in finite times. In the critical or subcritical case, we give a construction of the process coming down from infinity. Finally, it is shown that the continuous-state polynomial branching process arises naturally as the rescaled limit of a sequence of discrete-state processes.

\smallskip



\bigskip


\section{Introduction}

\setcounter{equation}{0}

Branching processes are models for the evolution of populations of particles. Those processes constitute an important subclass of Markov processes. Standard references on those processes with discrete-state space $\mbb{N}:= \{0,1,2,\dots\}$ are \citet{Har63} and \citet{AtN72}. As the quantity of particles can sometimes be expressed by other means than by counting, it is reasonable to consider branching models with continuous-state space $\mbb{R}_+:= [0,\infty)$. A diffusion process of that type was first studied by Feller (1951), which is now referred to as the \emph{Feller branching diffusion}. A general class of continuous-state branching processes were characterised in Lamperti (1967a) as the weak limits of rescaled discrete-state branching processes. The continuous-state branching models involve rich and deep mathematical structures and have attracted the attention of many researchers in the past decades. In particular, the connection of those processes with L\'evy processes through random time changes was pointed out by Lamperti (1967b). Multitype continuous-state branching processes were studied in Rhyzhov and Skorokhod (1970) and Watanabe (1969). A remarkable theory of flows of such processes with applications to flows of Bessel bridges and coalescents with multiple collisions has been developed by Bertoin and Le Gall (2000, 2003, 2005, 2006); see also Dawson and Li (2006, 2012). The reader may refer to Kyprianou (2006), Li (2011) and Pardoux (2016) for reviews of the literature in this subject.

It is well-known that the transition function $(P_t)_{t\geq 0}$ of a classical continuous-time branching process with state space $E=\mbb{N}$ or $\mbb{R}_+$ satisfies the following so-called \emph{branching property}:
 \beqlb\label{bp1}
P_t(x,\cdot)*P_t(y,\cdot)=P_t(x+y,\cdot),\qquad x,y\in E,
 \eeqlb
where ``$*$" denotes the convolution operation. The property means that different individuals in the population act independently of each other. In most realistic situations, however, this property is unlikely to be appropriate. In particular, when the number of particles becomes large or the particles move with high speed, the particles may interact and, as a result, the birth and death rates can either increase or decrease. Those considerations have motivated the study of \emph{generalised branching processes}, which may not satisfy (\ref{bp1}).

\subsection{Polynomial branching processes}

Let $\alpha$ and $b_i$, $i=0,1,\ldots$ be positive constants satisfying $b_1=0$ and $\sum^\infty_{i=0}b_i\le 1$. A \emph{discrete-state polynomial branching process} is a Markov chain on $\mbb{N}$ with $Q$-matrix $(q_{ij})$ defined by
 \beqlb\label{1.2aa}
q_{ij}=\left\{
\begin{array}{lcl}
 \alpha i^\theta b_{j-i+1}, & & {j\geq i+1, i\ge 1,}\\
 -\alpha i^\theta, & & {j=i\ge 1,}\\
 \alpha i^\theta b_0, & & {j=i-1, i\ge 1,}\\
 0, & & \mbox{otherwise.}
\end{array} \right.
 \eeqlb
Observe that $q_{ij}= i^\theta \rho_{ij}$, where $(\rho_{ij})$ is the $Q$-matrix of a random walk on the space of integers with jumps larger than $-1$. The transition rate of the discrete-state polynomial branching process is given by the power function $i\mapsto i^\theta$ and its transition distribution is given by the sequence $\{b_i: i\ge 0\}$. The $Q$-matrix (\ref{1.2aa}) is essentially a particular form of the model introduced by Chen (1997), who considered more general branching structures. Those processes have attracted the research interest of many other authors; see, e.g. Chen (2002), Chen et al.\ (2008) and Pakes (2007). When $\theta=1$, the model reduces to a classical discrete-state branching process, which satisfies property (\ref{bp1}). We refer to Chen (2004) for the general theory of continuous-time Markov chains.

In this paper, we introduce and study a continuous-state version of the process defined by (\ref{1.2aa}). Let $C^2_0[0,\infty)$ be the space of twice continuously differentiable functions on $[0,\infty)$ which together with their derivatives up to the second order vanish at $\infty$. By convention, we extend each $f$ on $[0,\infty)$ to $[0,\infty]$ by setting $f(\infty) = 0$. Fix a constant $\theta> 0$ and let
 \beqnn
\mcr{D}(L) = \Big\{f\in C^2_0[0,\infty): \lim_{x\to \infty} x^\theta |f^{(n)}(x)| = 0, n=0, 1, 2\Big\},
 \eeqnn
 where $f^{(n)}$ denotes the $n$-th derivative of $f$. Let $b\in\mbb{R}$ and $c\geq 0$ be constants and $m(du)$ a $\sigma$-finite measure on $(0,\infty]$ satisfying
 \beqnn
\int_{(0,\infty)} (1\wedge u^2) m(du)< \infty \qquad\mbox{and}\qquad m(\{\infty\})=a\geq 0.
 \eeqnn
For $x\in [0,\infty)$ and $f\in \mcr{D}(L)$ we define
 \beqlb\label{0.2}
Lf(x) = x^\theta \bigg[-af(x) - bf'(x) + cf''(x) + \int_{(0,\infty)} D_zf(x) m(dz)\bigg],
 \eeqlb
where
 \beqnn
D_zf(x) = f(x+z)-f(x)-zf'(x)1_{\{z\le 1\}}.
 \eeqnn
By Taylor's expansion one can see $Lf(\infty) := \lim_{x\to \infty}Lf(x) = 0$. In this work, a stochastically continuous Markov process $(X_t: t\ge 0)$ with state space $[0, \infty]$ is called a \emph{continuous-state polynomial branching process} if it has traps $0$ and $\infty$ and its transition semigroup $(P_t)_{t\geq 0}$ satisfies the forward Kolmogorov equation
 \beqlb\label{d1}
\frac{dP_tf}{dt}(x) = P_tLf(x),\qquad x\in [0,\infty], f\in \mcr{D}(L).
 \eeqlb

We call $\theta>0$ the \emph{rate power} of the continuous-state polynomial branching process. The ordinary continuous-state branching process corresponds to the special case $\theta=1$, which we refer to as the \emph{classical branching} case; see, e.g., Lamperti (1967a, 1967b). That is the only situation where the branching property (\ref{bp1}) is satisfied. Let $\psi$ be the function on $[0,\infty)$ defined by
 \beqlb\label{0.1}
\psi(\lambda) = -a+b\lambda+c\lambda^2+\int_{(0,\infty)}(e^{-\lambda z} - 1 + \lambda z1_{\{z\le 1\}}) m(dz),\qquad \lambda\ge 0.
 \eeqlb
We call $\psi$ the \emph{reproduction mechanism} of the process.  It can be easily checked that $\psi(0)= -a\leq 0$ and $\psi$ is a convex function. By (\ref{0.1}), dominated convergence and monotone convergence we see
\beqlb\label{0.2'}
\beta := \psi'(0) = b - \int_{(1,\infty)} z m(dz).
\eeqlb
Note that $-\psi$ is the Laplace exponent of a spectrally positive L\'{e}vy process. In this paper, we always assume that there exists $\lambda\in (0,\infty)$ such that $\psi(\lambda)> 0$, i.e., $-\psi$ does not corresponds to a subordinator.

\subsection{Construction of the process}

We here present a construction of the continuous-state polynomial branching process in terms of a stochastic equation with jumps. Suppose that $(\Omega,\mcr{F},\mcr{F}_t,\mbf{P})$ is a filtered probability space satisfying the usual hypotheses. Let $(B_t: t\ge 0)$ be an $(\mcr{F}_t)$-Brownian motion. Let $M(ds,dz,du)$ be an $(\mcr{F}_t)$-Poisson random measure on $(0,\infty)\times (0,\infty]\times (0,\infty)$ with intensity $dsm(dz)du$ and $\tilde{M}(ds,dz,du)$ the compensated measure. Let $X_0$ be a positive $\mcr{F}_0$-measurable random variable. We consider positive solutions of the stochastic integral equation
 \beqlb\label{e0.1}
X_t\ar=\ar X_0 +\sqrt{2c}\int^{t}_0 X_s^{\theta/2}dB_s + \int_0^t\int_{(0,1]} \int^{X_{s-}^\theta}_0 z \tilde{M}(ds,dz,du)\cr
 \ar\ar\quad
-\, b \int^{t}_0 X_s^\theta ds + \int_0^{t}\int_{(1,\infty]}\int^{X_{s-}^\theta}_0 z M(ds,dz,du).
 \eeqlb
Here and in the sequel, we understand $\int_a^b=\int_{(a,b]}$ and $\int_a^\infty=\int_{(a,\infty)}$ for $0\le a\le b< \infty$. For $\theta= 1$ the above equation has been considered in Dawson and Li (2006) and Fu and Li (2010). By saying $X= (X_t: t\geq 0)$ is a solution of the (\ref{e0.1}) we mean it is a c\`adl\`ag $[0,\infty]$-valued $(\mcr{F}_t)$-adapted process satisfying~(\ref{e0.1}) up to time $\zeta_n := \inf\{t\geq0: X_t\ge n$ or $X_t\le 1/n\}$ for each $n\ge 1$ and $X_t = \lim_{n\to \infty}X_{\zeta_n-}$ for $t\ge \tau:= \lim_{n\to \infty} \zeta_n$. Then both of the boundary points $0$ and $\infty$ are absorbing for $X$.

\btheorem\label{t1.1} {\rm(1)}~For any initial value $X_0=x\in [0,\infty)$ there exists a pathwise unique solution to {\rm(\ref{e0.1})}. {\rm(2)}~Let $X^x:= (X_t^x: t\ge 0)$ be the solution to {\rm(\ref{e0.1})} with $X_0=x$. Then $y\ge x\in [0,\infty)$ implies $\mbf{P}(X_t^y\ge X_t^x$ for every $t\ge 0) = 1$. {\rm(3)}~The solution to {\rm(\ref{e0.1})} is a continuous-state polynomial branching process defined by (\ref{d1}).
 \etheorem

 The stopping time $\tau$ is referred to as the \emph{absorbing time} of the continuous-state polynomial branching process. From (\ref{e0.1}) it is clear that $X$ does not have negative jumps, which is a crucial property. The process can also be obtained from a spectrally positive L\'{e}vy process through Lamperti type transformations, which play an important role in the study.

Let $D$ be the space of c\`{a}dl\`{a}g functions $w: [0,\infty)\to [0,\infty]$ with $0$ and $\infty$ as traps. Let $\rho(x,y)= |e^{-x}-e^{-y}|$ for all $x,y\in [0,\infty)$. We extend $\rho$ to a metric on $[0,\infty]$ making it homeomorphic to $[0,1]$. It is easy to see that $\rho(x,y)\le 1\land |x-y|$ for all $x,y\in [0,\infty)$. Then we define the uniform distance $\rho_\infty$ on $D$ by
 \beqlb\label{de0.1ss}
\rho_\infty(v,w) = \sup_{s\in[0,\infty)}\rho(v(s),w(s)),\qquad v,w\in D.
 \eeqlb
Let $\Lambda$ be the set of increasing homeomorphisms of $[0,\infty)$ into itself and define the metric $d_\infty$ on $D$ by
 \beqlb\label{de0.1}
d_\infty(v,w) := \inf_{\lambda\in\Lambda}\rho_\infty(v,w\circ\lambda) \vee \|\lambda-I\|,\qquad v,w\in D,
 \eeqlb
where $I$ is the identity and $\|\cdot\|$ is the supremum norm. For $x\in [0,\infty)$ let $P^x$ denote the distribution on $D$ of the process $X^x= (X^x_t: t\ge 0)$ defined by (\ref{e0.1}) with initial value $X_0=x$.

\btheorem\label{t1.4xx} The mapping $[0,\infty)\ni x\mapsto P^x$ is continuous in the weak convergence topology. \etheorem

For any $y\in [0,\infty]$ let $\tau_y= \inf\{t\geq 0: X_t=y\}$. We call $\tau_0$ the \emph{extinction time} and $\tau_\infty$ the \emph{explosion time} of the continuous-state polynomial branching process $X$. Note that $\tau_0\wedge \tau_\infty= \tau$ is the absorbing time. For $x\in [0,\infty]$ write $\mbf{P}_x= \mbf{P}(\,\cdot\,|X_0=x)$. Since $X$ has no negative jump, we have $\mbf{P}_x(X_{\tau_y}= y)=1$ for $x\ge y\in [0,\infty)$. Let $q= \inf\{\lambda>0: \psi(\lambda)> 0\}$ be the largest root of $\psi(\lambda)=0$.

\btheorem\label{t1.4} {\rm(1)} For any $y\le x\in (0,\infty)$ we have
 \beqlb\label{t1.4sss}
\mbf{P}_x(\tau_y< \infty)= e^{-q(x-y)}.
 \eeqlb
{\rm(2)} For any $x\in [0,\infty)$ we have
 \beqlb\label{2.10}
\mbf{P}_x(X_\infty = 0) = e^{-qx},
 \quad
\mbf{P}_x(X_\infty= \infty)= 1- e^{-qx},
 \eeqlb
where $X_\infty=\lim_{t\to\infty}X_t.$
 \etheorem

By applying Theorem~7.8 in Either and Kurtz (1986, p. 131) and Theorem~\ref{t1.4xx} we see the process is a Feller process. From the transition semigroup $(P_t)_{t\ge 0}$ of the process we define its \emph{resolvent} $(U^\eta)_{\eta>0}$ by
 \beqlb\label{1.16}
U^\eta(x,dy)= \int_0^\infty e^{-\eta t}P_t(x,dy)dt, \qquad x,y\in [0,\infty].
 \eeqlb
The next theorem gives a characterization of the resolvent and plays the key role in the study of the hitting times of $X$. Let $e_\lambda(x)= e^{-\lambda x}$ for $\lambda\in (0,\infty)$ and $x\in [0,\infty]$.

\btheorem\label{t4.1} For any $\eta,\lambda>0$ and $x\in [0,\infty)$ we have
 \beqlb\label{2.11}
\eta U^\eta e_\lambda(x) - e^{-\lambda x}
 =
\psi(\lambda)\int_{[0,\infty)} y^\theta e^{-\lambda y} U^\eta(x,dy)
 \eeqlb
and
 \beqlb\label{1.11xx}
\int^\infty_\lambda l_x(\eta,z)(z-\lambda)^{\theta-1} dz
=
\Gamma(\theta)\int^\infty_0 e^{-\eta t}dt\int_{(0,\infty)}e^{-\lambda y} P_t(x,dy),
 \eeqlb
 where $\Gamma$ denotes the Gamma function and, for $\eta, z>0$,
\beqlb\label{de0.2}
l_x(\eta,z) = \int_{[0,\infty)} y^\theta e^{-zy} U^\eta(x,dy)
 =
\psi(z)^{-1}\big[\eta U^\eta e_z(x)-e^{-zx}\big].
 \eeqlb
\etheorem

We remark that, in the case $q>0$, the second expression in (\ref{de0.2}) should be understood by continuity at $z=q$.

\subsection{Mean extinction and explosion times}

In this subsection, we give some expressions for the mean hitting times of the continuous-state polynomial branching process $X= (X_t: t\ge 0)$. Let $\mbf{E}_x$ denote the expectation with respect to the conditional law $\mbf{P}_x= \mbf{P}(\,\cdot\,|X_0=x)$. Recall that $q= \inf\{\lambda>0: \psi(\lambda)> 0\}$ is the largest root of $\psi(\lambda)=0$. Let
 \beqlb\label{de0.3}
h_x(\lambda)= \frac{e^{-qx} - e^{-\lambda x}}{\psi(\lambda)}, \qquad \lambda>0, x>0
 \eeqlb
with $h_x(q)= xe^{-qx}/\psi'(q)$ by continuity if $q>0$.

 \btheorem\label{t1.7}
For any $x\in (0,\infty)$ we have the moment formulas:
 \beqlb\label{2.3}
\mbf{E}_x(\tau_0: X_\infty= 0)
 =
\frac{1}{\Gamma(\theta)} \int^\infty_0 h_x(\lambda+q)\lambda^{\theta-1}d\lambda,
 \eeqlb
 \beqlb\label{2.3'}
\mbf{E}_x(\tau_\infty: X_\infty= \infty)
 =
\frac{1}{\Gamma(\theta)} \int^\infty_0[h_x(\lambda) - h_x(\lambda+q)] \lambda^{\theta-1}d\lambda.
 \eeqlb
and
 \beqlb\label{2.4}
\mbf{E}_x(\tau)= \mbf{E}_x(\tau_\infty\wedge\tau_0) = \frac{1}{\Gamma(\theta)}\int^\infty_0 h_x(\lambda)\lambda^{\theta-1}d\lambda.
 \eeqlb
\etheorem

\btheorem\label{t1.8} For any $y\le x\in (0,\infty)$ we have
 \beqlb\label{0.1''}
\mbf{E}_x(\tau_\infty\wedge\tau_y)
 =
\frac{\,e^{-qx}\,}{\Gamma(\theta)}\int_0^\infty\frac{e^{-(\lambda-q)y} - e^{-(\lambda-q)x}}{\psi(\lambda)}\lambda^{\theta-1}d\lambda.
 \eeqlb
\etheorem

\bcorollary\label{t1.6'} Suppose that $a= -\psi(0)= 0$ and $\psi'(0)\geq 0$. Then for $y\le x\in (0,\infty)$ we have
 \beqlb\label{2.9}
\mbf{E}_x(\tau_y) = \frac{1}{\Gamma(\theta)}\int^\infty_0 \frac{e^{-\lambda y}-e^{-\lambda x}}{\psi(\lambda)} \lambda^{\theta-1} d\lambda.
 \eeqlb
 \ecorollary

The discrete-state versions of (\ref{2.3}) and (\ref{2.3'}) were proved in Chen (2002) and Pakes (2007), respectively. As far as we know, the discrete-state form of (\ref{0.1''}) has not been established in the literature. One may compare (\ref{2.9}) with Corollary~9 in Duhalde et al.\ (2014). It seems other moment formulas are new also for classical continuous-state branching processes.

\subsection{Extinction and explosion probabilities}

The two theorems presented in this subsection are about the extinction and explosion probabilities of the process. They generalize the results in Grey (1974) and Kawazu and Watanabe (1971), where the classical branching case $\theta=1$ was studied.

\btheorem\label{t1.4'} {\rm(1)} In the case $\theta\geq 2$, for any $x\in (0,\infty)$ we have $\mbf{P}_x(\tau_0< \infty)= 0$.
{\rm(2)} In the case $0<\theta< 2$, for any $x\in (0,\infty)$ we have $\mbf{P}_x(\tau_0< \infty)> 0$ if and only if
 \beqlb\label{1.1'}
\int^\infty_\varepsilon \frac{\lambda^{\theta-1}}{\psi(\lambda)} d\lambda <\infty
 \eeqlb
for some and hence all $\varepsilon\in (q,\infty)$. {\rm(3)} If $\mbf{P}_x(\tau_0< \infty)> 0$, then $\mbf{P}_x(\tau_0< \infty)= \mbf{P}_x(X_\infty= 0)= e^{-qx}.$
\etheorem

\bcorollary\label{r1.4'} If $0<\theta<2$ and $c>0$, for any $x\in (0,\infty)$ we have $\mbf{P}_x(\tau_0< \infty)> 0$. \ecorollary

\btheorem\label{t1.3} {\rm(1)} In the case $\theta> 1$, for any $x\in (0,\infty)$ we have $\mbf{P}_x(\tau_\infty< \infty)= 0$ if and only if $a= -\psi(0)= 0$ and $\psi'(0)\geq 0$.
{\rm(2)} In the case $0<\theta\le 1$, for any $x\in (0,\infty)$ we have $\mbf{P}_x(\tau_\infty< \infty)= 0$ if and only if $a= -\psi(0)= 0$ and one of the following two conditions is satisfied: {\rm(i)} $\psi'(0)>-\infty$; {\rm(ii)} $\psi'(0)= -\infty$ and
\beqlb\label{1.1}
\int_0^\varepsilon \frac{\lambda^{\theta-1}}{-\psi(\lambda)} d\lambda= \infty.
 \eeqlb
for some and hence all $\varepsilon\in (0,q)$.
\etheorem

\subsection{The process coming down from infinity}

Let $X^x= (X^x_t: t\geq 0)$ be defined as in Theorem~\ref{t1.1}. Since the process has no negative jump, for any $x\ge y\in [0,\infty)$ we can define $\tau^x_y= \inf\{t\geq 0: X^x_t= y\}$. By Theorem~\ref{t1.1}\,(2), we see the mapping $x\mapsto \tau^x_y$ is increasing in $x\in [y,\infty)$, thus the limit $\tau^\infty_y:= \lim_{x\rightarrow \infty}\tau^x_y$ exists. It is easy to see that $y\mapsto \tau^\infty_y$ is decreasing in $y\in [0,\infty)$. By Corollary~\ref{t1.6'} for any $y\in (0,\infty)$ we have
 \beqlb\label{1.5}
\mbf{E}(\tau^\infty_y)
 =
\lim_{x\to\infty} \mbf{E}(\tau^x_y)
 =
\frac{1}{\Gamma(\theta)}\int^\infty_0 \frac{e^{-\lambda y}}{\psi(\lambda)} \lambda^{\theta-1} d\lambda.
 \eeqlb

\btheorem\label{t3.2} Suppose that $a= -\psi(0)= 0$ and $\psi'(0)\geq 0$. Then the following four statements are equivalent:
 \benumerate

\itm[\rm(i)] $\mbf{P}(\tau^{\infty}_y< \infty)> 0$ for each $y\in (0,\infty)$;

\itm[\rm(ii)] $\mbf{P}(0<\tau^{\infty}_y< \infty)= 1$ for each $y\in (0,\infty)$;

\itm[\rm(iii)] $\mbf{E}(\tau^{\infty}_y)< \infty$ for each $y\in (0,\infty)$;

\itm[\rm(iv)] for each $\varepsilon\in (0,\infty)$ we have
 \beqlb\label{1.4}
\int_0^\varepsilon \frac{\lambda^{\theta-1}}{\psi(\lambda)} d\lambda< \infty.
 \eeqlb
 \eenumerate
\etheorem

\bcorollary\label{c3} Suppose that $a= -\psi(0)= 0$, $\psi'(0)\geq 0$ and (\ref{1.4}) holds. Then we have a.s.\ $\lim_{y\to\infty} \tau^\infty_y= 0$
\ecorollary

By saying a process $(X_t)_{t> 0}$ is a solution to (\ref{e0.1}) with initial state $\infty$, we mean $\lim_{t\downarrow 0} X_t=\infty$ and, for $t> r> 0$,
\beqlb\label{ne1}
X_t\ar=\ar X_r +\sqrt{2c}\int^{t}_r\sqrt{ X_s^{\theta}}dB_s + \int_r^t\int_{(0,1]} \int^{{X}_{s-}^\theta}_0 z \tilde{M}(ds,dz,du)\cr
 \ar\ar\quad
-\, b \int^{t}_r X_s^\theta ds + \int_r^{t}\int_{(1,\infty]}\int^{X_{s-}^\theta}_0 z M(ds,dz,du).
\eeqlb

\btheorem\label{t3.3}
Suppose that $a= -\psi(0)= 0$, $\psi'(0)\geq 0$ and (\ref{1.4}) holds. Then: {\rm(i)} there is a pathwise unique solution $(X^\infty_t)_{t> 0}$ to {\rm(\ref{e0.1})} with initial state $\infty$; {\rm(ii)} $\tau^\infty_y = \inf\{t>0: X^\infty_t=y\}$ for $y\in (0,\infty)$; {\rm(iii)} if we set $X^\infty_0= \infty$, then $(X^x_t)_{t\ge 0}$ converges a.s.\ to $(X_t^\infty)_{t\ge 0}$ in $(D,\rho_\infty)$ as $x\to \infty$.
\etheorem

The above theorem shows that a solution of (\ref{e0.1}) may come down from $\infty$. This property is not possessed by classical branching processes. In fact, to guarantee the integrability (\ref{1.4}) we should at least have $\theta>1$. We refer to Bansaye et al.\ (2015) for a study of the speed of coming down from infinity of birth and death process. For coalescent processes and branching models with interaction, similar phenomena have been observed and studied by a number of authors; see, e.g., Berestycki et al.\ (2010, 2014), Lambert (2005) and Pardoux (2016) and the references therein.

\subsection{Convergence of discrete-state processes}

The following theorem shows that the continuous-state polynomial branching process $X = (X_t: t\ge 0)$ defined by (\ref{e0.1}) can be obtained as the limit of a sequence of rescaled discrete-state branching processes.

\btheorem\label{t3.1}
There exists a sequence of discrete-state polynomial branching processes $\xi_n= (\xi_n(t): t\ge 0)$ and a sequence of positive number $\gamma_n$, $n= 1,2,\dots$ such that $(n^{-1}\xi_n(\gamma_nt): t\ge 0)$ converges to $(X_t: t\ge 0)$ weakly in $(D,d_\infty)$.
\etheorem

\bigskip

To conclude the introduction, we give the following two examples of the continuous-state polynomial branching process.

\bexample\label{t3.1zz} Let $0<\theta<1$ and consider the equation
 \beqlb\label{ex1}
X_t = \int^t_0 X^\theta_t dt, \qquad z(0)=0.
 \eeqlb
Obviously $X_t=(1-\theta)^{1/(1-\theta)}t^{1/1-\theta}$ is a solution to~(\ref{ex1}). It is trivial to see that $x_0(t)\equiv 0$ is another solution to the above equation. Then the requirement of $0$ being a trap is necessary to guarantee the pathwise uniqueness of the solution to~(\ref{e0.1}). \qed \eexample

\bexample\label{t3.1yy} A continuous-state polynomial branching process with reproduction mechanism $\psi(\lambda)= c\lambda^2 + b\lambda$ ($c>0$, $b>0$) is defined by
 \beqlb\label{e0.1ab}
X_t\ar=\ar X_0 +\sqrt{2c}\int^t_0 X_s^{\theta/2}dB_s - b\int^t_0 {X_s^{\theta}}ds.
 \eeqlb
For this process we have a.s.\ $\tau_\infty=\infty$ and the formulas given above take simple forms. From (\ref{1.5}) we have, for $y\in (0,\infty)$,
 \beqnn
\mbf{E}(\tau^\infty_y)
 =
\frac{1}{\Gamma(\theta)}\int^\infty_0 \frac{e^{-\lambda y}\lambda^{\theta-1}} {c\lambda^2 + b\lambda} d\lambda,
 \eeqnn
which is finite if and only if $\theta>1$. By letting $y\to 0$ in the above equality we get
 \beqnn
\mbf{E}(\tau^\infty_0)
 =
\frac{1}{\Gamma(\theta)}\int^\infty_0 \frac{\lambda^{\theta-1}} {c\lambda^2 + b\lambda} d\lambda,
 \eeqnn
which is finite if and only if $1<\theta<2$. The above formula gives explicitly the expected time for the process $X$ to cross the state space from $\infty$ to $0$. The process defined by (\ref{e0.1ab}) reduces to a classical Feller branching diffusion when $\theta=1$. A closely related model has been studied recently by Berestycki et al.\ (2015). The polynomial branching structure has also appeared in the so-called generalized Cox-Ingersoll-Ross model; see, for example, Borkovec and Kl\"uppelberg (1998) and Fasen et al.\ (2006). \qed \eexample

We present the proofs of the results in the following sections. Section~2 is devoted to the construction of the process. The mean extinction and explosion times are  computed in Section~3. In Section~4, the extinction and explosion probabilities are explored. In Section~5, we prove the construction of the process coming down from $\infty$. The convergence of discrete-state processes is discussed in Section~6.

\section{Construction of the process}

\setcounter{equation}{0}

In this section, we construct the continuous-state polynomial branching process $X$ in terms of stochastic equations and random time changes.

\medskip

\noindent\emph{Proof of Theorem}~\ref{t1.1}.~ (1) We prove the result by an approximation argument. For each $n\ge 1$ define
 \beqlb\label{0.4}
r_n(x)=\left\{
\begin{array}{lcl}
 n^{\theta}, & & n< |x|< \infty, \cr
 |x|^\theta, & & 1/n< |x|\leq n, \cr
 n^{2-\theta}|x|^2, & & 0\le |x|\leq 1/n.
\end{array} \right.
 \eeqlb
By Theorem 9.1 in Ikeda and Watanabe (1989, p.245) there is a pathwise unique solution $\{\xi_n(t): t\ge 0\}$ to the stochastic equation
 \beqlb\label{e0.2}
\xi_n(t) \ar=\ar x + \int^t_0 \sqrt{2cr_n(\xi_n(s))} dB(s) + \int^t_0 \int_{(0,1]} \int^{r_n(\xi_n(s-))}_0 z \tilde{M}(ds,dz,du)\cr
 \ar\ar\quad
- \int^t_0 br_n(\xi_n(s))ds + \int^t_0 \int_{(1,\infty]}\int^{r_n(\xi_n(s-))}_0 (z\land n) M(ds,dz,du).
 \eeqlb
Let $\zeta_n= \inf\{t\ge 0: \xi_n(t)\geq n$ or $\xi_n(t)\leq 1/n\}$. Clearly, the sequence of stopping times $\{\zeta_n\}$ is increasing and $\xi_n(t) = \xi_m(t)$ for $t\in [0,\zeta_{m\wedge n})$. Let $\tau = \lim_{n\to \infty} \zeta_n$. We define the process $(X_t: t\ge 0)$ by $X_t = \xi_n(t)$ for $t\in [0,\zeta_n)$ and $X_t = \lim_{n\to \infty}\xi_n(\zeta_n)$ for $t\in [\tau,\infty)$. Then $\zeta_n= \inf\{t\ge 0: X_t\geq n$ or $X_t\leq 1/n\}$ and $(X_t: t\ge 0)$ is a solution of (\ref{e0.1}). The pathwise uniqueness of the solution follows from that for (\ref{e0.2}) in the time interval $[0,\zeta_n)$ for each $n\ge 1$.

(2) Let $\{\xi_n^x(t): t\geq 0\}$ denote the solution of (\ref{e0.2}) to indicate its dependence on the initial state. For any $y\ge x\ge 0$, we can use Theorem~5.5 in Fu and Li (2010) to see $P(\xi_n^y(t)\ge \xi_n^x(t)$ for every $t\ge 0)=1$, and so $P(X_t^y\ge X_t^x$ for every $t\ge 0)= 1$.

(3) For any $t\ge 0$ let $P_t(x,\cdot)$ be the distribution of $X_t$ on $[0,\infty]$ with $X_0= x\in [0,\infty]$. By Theorem~\ref{t1.1}, for any $y\in [0,\infty]$, the mapping $x\mapsto P_t(x,[0,y])$ is decreasing, so it is Borel measurable. A monotone class argument shows $x\mapsto P_t(x,A)$ is Borel measurable for each Borel set $A\subset [0,\infty]$. Then $P_t(x,dy)$ is a Borel kernel on $[0,\infty]$. For any finite $(\mcr{F}_t)$-stopping time $\sigma$, from the equation (\ref{e0.1}) we have
 \beqnn
X_{\sigma+t}\ar=\ar X_\sigma +\sqrt{2c}\int^{t}_0 X_{\sigma+s-}^{\theta/2}dB_{\sigma+s} + \int_0^t\int_{(0,1]} \int^{X_{\sigma+s-}^\theta}_0 z \tilde{M}(\sigma+ds,dz,du)\cr
 \ar\ar\qquad
-\, b \int^t_0 X_{\sigma+s-}^\theta ds + \int_0^t\int_{(1,\infty]}\int^{X_{\sigma+s-}^\theta}_0 z M(\sigma+ds,dz,du).
 \eeqnn
Here $(B_{\sigma+s}-B_\sigma: s\ge 0)$ is a Brownian motion and $M(\sigma+ds,dz,du)$ is a Poisson random measure on $(0,\infty)^3$ with intensity $dsm(dz)du$. Those are true under the original probability $\mbf{P}(\cdot)$ and also under the conditional probability $\mbf{P}(\cdot|\mcr{F}_\sigma)$ because of the independent increment property. In particular, under $\mbf{P}(\cdot|\mcr{F}_\sigma)$, the process $(X_{\sigma+t}: t\ge 0)$ satisfies a stochastic equation of the same form as (\ref{e0.1}) with initial state $X_\sigma$. From the uniqueness of solution of (\ref{e0.1}), it follows that $\mbf{P}(X_{\sigma+t}\in \cdot|\mcr{F}_\sigma)= P_t(X_\sigma,\cdot)$. That gives the strong Markov property of the process $(X_t: t\ge 0)$. For $f\in \mcr{D}(L)$, we can use (\ref{e0.2}) and It\^o's formula to see
 \beqnn
f(\xi_n(t\wedge \zeta_n))
 \ar=\ar
f(x) + b\int^{t\wedge \zeta_n}_0 f'(\xi_n(s))r_n(\xi_n(s))ds + c\int^{t\wedge \zeta_n}_0 f''(\xi_n(s))r_n(\xi_n(s)) ds \cr
 \ar\ar
+ \int^{t\wedge \zeta_n}_0 r_n(\xi_n(s))ds\int_{(0,\infty]} \Big[f(\xi_n(s)+z) - f(\xi_n(s)) \cr
 \ar\ar\qqquad\qqquad\qqquad
-\, f'(\xi_n(s))z 1_{\{z\le 1\}}\Big] m(dz) + M_n(t) \cr
 \ar=\ar
f(x) + \int^{t\wedge \zeta_n}_0 Lf(\xi_n(s))ds + M_n(t),
 \eeqnn
where
\beqnn
 M_n(t)\ar=\ar\int^{t\wedge \zeta_n}_0 f'(\xi_n(s))\sqrt{2cr_n(\xi_n(s))} dB(s)\cr
 \ar\ar
+ \int^{t\wedge \zeta_n}_0 \int^{r_n(\xi_n(s))}_0\int_{(0,\infty]} [f(\xi_n(s)+z) - f(\xi_n(s))] \tilde{M}(ds,dz,du).
\eeqnn
Notice that $(M^{(n)}_f(t): t\ge0)$ is a martingale bounded on each bounded time interval. Then letting $n\to \infty$ in the above equality and using bounded convergence theorem we get
 \beqnn
f(X_{t\wedge \tau})= f(x) + \int^{t\wedge \tau}_0 Lf(X_s)ds + \mbox{martingale}.
 \eeqnn
Since $0$ and $\infty$ are traps for $(X_t: t\ge 0)$ and $Lf(0) = Lf(\infty) = 0$, it follows that
 \beqlb\label{2.6}
f(X_t)= f(x) + \int^t_0 Lf(X_s)ds + \mbox{martingale}.
 \eeqlb
Then we take the expectation on both sides and obtain
 \beqnn
P_tf(x)= f(x) + \int^t_0 P_sLf(x)ds.
 \eeqnn
That gives the Kolmogorov forward equation (\ref{d1}). \qed

We next give some results connecting the process $X$ and a spectrally positive L\'{e}vy process through Lamperti type transformations. This connection plays an important role in the study of properties of the process. Let $Z= (Z_t: t\ge 0)$ be a spectrally positive L\'{e}vy process with \emph{Laplace exponent} $-\psi$ and initial state $Z_0=x\ge 0$. Note that $Z$ is absorbed by $\infty$ after an exponential time $T_\infty$ with parameter $a= -\psi(0)\geq 0$. Let $T_y= \inf\{t\geq 0: Z_t=y\}$ for $y\in [0,\infty]$. Let $T = T_0\wedge T_\infty$ be the \emph{absorbing time} of $Z$. Let $Y_t= Z_{t\wedge T}$ for $t\ge 0$. We call $Y:= (Y_t: t\ge 0)$ an \emph{absorbed spectrally positive L\'{e}vy process}. By Proposition 37.10 in Sato (1999, p.255), the limit $Y_\infty:= \lim_{t\to \infty} Y_t$ exists a.s.\ in $[0,\infty]$.

\bproposition\label{t1.0}
Let $\alpha(t)= \int^t_0 Y^{-\theta}_{s-}ds$ and $\eta(t)= \inf\{s\geq 0: \alpha(s)> t\}$ for $t\ge 0$. Then $J_\theta(Y):= (Y_{\eta(t)}: t\ge 0)$ solves {\rm(\ref{e0.1})} on an extension of the original probability space. \eproposition

\proof Let $W(t)$ be a Brownian motion and let $N_0(ds,dz)$ be a Poisson random measure on $(0,\infty)\times (0,\infty]$ with intensity $dsm(dz)$. Then a realization of the L\'{e}vy process $Z := (Z_t: t\ge 0)$ is defined by
 \beqlb\label{0.3}
Z_t = x - bt + \sqrt{2c} W(t) + \int^t_0 \int_{(0,1]} z \tilde{N}_0(ds,dz)
+\int^t_0 \int_{(1,\infty]} z N_0(ds,dz),
 \eeqlb
where $\tilde{N}_0(ds,dz)= N_0(ds,dz)-dsm(dz)$. Let $\{(s_i,z_i): i=1,2,\dots\}$ be an enumeration of the atoms of $N_0(ds,dz)$. On an extension of the original probability space, we can construct a sequence of $(0,1]$-valued i.i.d.\ uniform random variables $\{u_i\}$ independently of $W(t)$ and $N_0(ds,dz)$. Then
 $$
 M_0(ds,dz,du) := \sum^{\infty}_{i=1} \delta_{(s_i,z_i,u_i)}(ds,dz,du)
 $$
defines a Poisson random measure on $(0, \infty)\times(0,\infty]\times (0,1]$ with intensity $ds m(dz) du$. Let $\tilde{M}_0(ds,dz,du) = M_0(ds,dz,du)- dsm(dz)du$. Then we have
 \beqnn
Z_t \ar=\ar x - bt + \sqrt{2c}W(t) + \int^t_0\int_{(0,1]}\int^1_0 z \tilde{M}_0(ds,dz,du)\cr
 \ar\ar\qquad
+ \int^t_0\int_{(1,\infty]}\int^1_0 z M_0(ds,dz,du).
 \eeqnn
Let $Y= (Y_t: t\ge 0)$ be the absorbed process associated with $Z$ and let $X_t= Y_{\eta(t)}$ for $t\ge 0$. Let $\zeta_n= \inf\{t\ge 0: X_t\geq n$ or $X_t\leq 1/n\}$ and $\zeta = \lim_{n\to \infty} \zeta_n$. Then we have
 \beqlb\label{1.2}
X_{t\land \zeta_n} \ar=\ar x + \sqrt{2c} W(\eta(t\land \zeta_n)) + \int^{\eta(t\land \zeta_n)}_0 \int_{(0,1]} \int_0^1 z \tilde{M}_0(ds,dz,du) - b\eta(t\land \zeta_n) \cr
 \ar\ar\qquad
+ \int^{\eta(t\land \zeta_n)}_0 \int_{(1,\infty]} \int^1_0 z M_0(ds,dz,du)\cr
 \ar=\ar
x + \sqrt{2c} W(\eta(t\land \zeta_n)) + \int^t_0 \int_{(0,1]} \int_0^1 z \tilde{M}_0(d\eta(s\land \zeta_n),dz,du) \cr
 \ar\ar\qquad
-\, b\eta(t\land \zeta_n) + \int^t_0 \int_{(1,\infty]} \int^1_0 z M_0(d\eta(s\land \zeta_n),dz,du).
 \eeqlb
By the definition of $\alpha(t)$ we have $d\alpha(t)= Y^{-\theta}_{t-}dt$ for $0\le t< T$ and $d\eta(t)= Y^\theta_{\eta(t)-}dt= X^\theta_{t-}dt$ for $0\le t< \tau$, where $\tau$ is defined after (\ref{e0.1}). It follows that
 \beqlb\label{Lam1.0}
\eta(t\land \zeta_n)
 =
\int^{t\land \zeta_n}_0 X_{s-}^\theta 1_{\{s< \tau\}}ds
 =
\int^{t\land \zeta_n}_0 X_{s-}^\theta ds.
 \eeqlb
By representation of time-changed Brownian motions, there is a Brownian motion $\{B(t)\}$ on an extension of the original probability space so that
 $$
W(\eta(t\land \zeta_n)) = \int^{t\land \zeta_n}_0 X^{\theta/2}_{s-}dB(s);
 $$
 see e.g. Theorem 4.3 in Ikeda and Watanabe (1989, p.198). On the extended probability space, we can take another independent Poisson random measure $\{ M_1(ds,dz,du)\}$ on $(0,\infty)\times(0,\infty]\times(0,\infty)$ with intensity $dsm(dz)du$ and define the random measure
 $$
M(ds,dz,du) = 1_{\{s< \tau,u\leq X^{\theta}_{s-}\}} M_0(d\eta(s),dz, X_{s-}^{-\theta}du)+ 1_{\{u> X^{\theta}_{s-}\}} M_1(ds,dz,du).
 $$
Using~(\ref{Lam1.0}) one can see $\{M(ds,dz,du)\}$ has the deterministic compensator $dsm(dz)du$, so it is a Poisson random measure. Then~(\ref{e0.1}) follows by substituting $ W(t)$ and $M_0(ds,dz,du)$ in (\ref{1.2}) and taking $n\to \infty$. \qed

\bproposition\label{t1.0'}
Let $\gamma(t) = \int^t_0 X^{\theta}_{s-}ds$ and $\beta(t)= \inf\{s\geq 0: \gamma(s)> t\}$ for $t\ge 0$. Then $L_{\theta}(X):= (X_{\beta(t)}: t\ge 0)$ is an absorbed spectrally positive L\'{e}vy process. \eproposition

\proof Without loss of generality, we assume $X_0=x\in [0,\infty)$ is deterministic. Let $Y_t= X_{\beta(t)}= X_{\beta(t)\land \tau}$ for $t\ge 0$ and let $T= \inf\{t\ge 0: Y_t= 0$ or $Y_t= \infty\}$. We have
 \beqlb\label{1.22}
Y_t \ar=\ar x + \sqrt{2c}\int^{\beta(t)\land \tau}_0 X_{s-}^{\theta/2}dB(s) + \int_0^{\beta(t)\land \tau} \int_{(0,1]}\int^{X_{s-}^\theta}_0 z \tilde{M}(ds,dz,du) \cr
 \ar\ar\quad
-\, b \int^{\beta(t)\land \tau}_0 X_{s-}^\theta ds + \int_0^{\beta(t)\land \tau}\int_{(1,\infty]} \int^{X_{s-}^\theta}_0 z M(ds,dz,du)\cr
 \ar=\ar
x + \sqrt{2c}\int^t_0 X_{\beta(s)-}^{\theta/2}1_{\{\beta(s)\le \tau\}}dB(\beta(s)) + \int_0^t \int_{(0,1]}\int^{X_{\beta(s)-}^\theta}_0 z1_{\{\beta(s)\le \tau\}} \tilde{M}(d\beta(s),dz,du) \cr
 \ar\ar\quad
-\, b \int^t_0 X_{\beta(s)-}^\theta1_{\{\beta(s)\le \tau\}} d\beta(s) + \int_0^t\int_{(1,\infty]} \int^{X_{\beta(s)-}^\theta}_0 z1_{\{\beta(s)\le \tau\}} M(d\beta(s),dz,du)\cr
 \ar=\ar
x + \sqrt{2c}\int^t_0 Y^{\theta/2}_{s-}1_{\{\beta(s)\le \tau\}} dB(\beta(s)) + \int_0^{t}\int_{(0,1]}\int^{Y_{s-}^\theta}_0 z1_{\{\beta(s)\le \tau\}} \tilde{M}(d\beta(s),dz,du)\cr
 \ar\ar\quad
-\, b \int^{t}_0 Y_{s-}^\theta1_{\{\beta(s)\le \tau\}} d\beta(s) + \int_0^t\int_{(1,\infty]} \int^{Y_{s-}^\theta}_0 z1_{\{\beta(s)\le \tau\}} M(d\beta(s),dz,du).
 \eeqlb
By the definition of $\gamma(t)$ and $\beta(t)$ we have $d\gamma(t) = X^\theta_{t-}dt$ for $0\le t< \tau$ and $d\beta(t) = X^{-\theta}_{\beta(t)-}dt = Y^{-\theta}_{t-}dt$ for $0\le t< T$. Thus
 $$
\int^t_0 Y^\theta_{s-}1_{\{\beta(s)<\tau\}}d\beta(s)
 =
\int^t_0 1_{\{s<T\}} ds = t\wedge T.
 $$
It follows that
 $$
W_0(t) := \int^t_0 Y^{\theta/2}_{s-}1_{\{\beta(s)<\tau\}} dB(\beta(s))
 $$
defines a continuous local martingale with $\langle W_0\rangle (t) = t\wedge T$. Then we can extend $\{W_0(t)\}$ to a Brownian motion $\{W(t)\}$. Now define the random measure $\{N_0(ds,dz)\}$ on $(0,\infty)\times(0,\infty]$ by
 $$
N_0((0,t]\times (a_1,a_2])= \int^t_0\int^{a_2}_{a_1}\int_0^{Y^\theta_{s-}}1_{\{\beta(s)<\tau\}} M(d\beta(s),dz,du),
 $$
where $t\geq 0$ and $a_1, a_2\in(0,\infty]$. It is easy to check that $\{N_0(ds,dz)\}$ has predictable compensator $Y^\theta_{s-}1_{\{\beta(s)<\tau\}}d\beta(s)m(dz) = 1_{\{s< T\}} dsm(dz)$. Then we can extend $\{N_0(ds,dz)\}$ to a Poisson random measure $\{N(ds,dz)\}$ on $(0,\infty)^2$ with intensity $dsm(dz)$; see, e.g., Ikeda and Watanabe (1989, p.93). From~(\ref{1.22}) it follows that
 \beqnn
Y_t = x +\sqrt{2c}W(t\wedge T) + \int^{t\wedge T}_0\int_{(0,1]}z \tilde{N}(ds,dz) - b (t\wedge T) + \int^{t\wedge T}_0\int_{(1, \infty]}z N(ds,dz).
 \eeqnn
Then $(Y_t)$ is an absorbed spectrally positive L\'{e}vy process. \qed

We call $L_{\theta}$ a \textit{generalized Lamperti transformation} and $J_\theta$ the \textit{inverse generalized Lamperti transformation}. In the particular case $\theta=1$, they reduce to the classical transformations introduced by Lamperti (1967a, 1967b).

\medskip

\noindent\emph{Proof of Theorem}~\ref{t1.4xx}.~ Let $Z = (Z_t: t\ge 0)$ be the L\'{e}vy process starting at $0$ with Laplace exponent $-\psi$. Let $Z_t^x= x+Z_t$ for $x\in [0,\infty)$. Let $T_0^x= \inf\{t\geq 0: Z_t^x= 0\}$ and $T_\infty^x= \inf\{t\geq 0: Z_t^x= \infty\}$. Then $x\mapsto T_0^x$ is a.s.\ increasing. By  Corollary 3.13 in Kyprianou (2006, p.82), for any $\lambda> 0$ we have
 \beqnn
\mbf{E}(e^{-\lambda T_0^x}) = \mbf{E}(e^{-\lambda T_0^x}1_{\{T_0^x< T_\infty^x\}})
 =
\exp\{-\psi^{-1}(\lambda)x\},
 \eeqnn
where $\psi^{-1}(\lambda) = \inf\{z\ge 0: \psi(z)> \lambda\}$. Then $\lim_{y\to x}T_0^y= T_0^x$ first in distribution and then almost surely. Let $Y^x = (Z_{t\wedge T_0^x}^x: t\ge 0)$. For $x< y\in [0,\infty)$ we have
 \beqnn
\rho(Z_{t\wedge T_0^x}^x, Z_{t\wedge T_0^y}^y)
 \ar\leq\ar
\rho(Z_{t\wedge T_0^x}^x, Z_{t\wedge T_0^x}^y) + \rho(Z_{t\wedge T_0^x}^y, Z_{t\wedge T_0^y}^y) \\
 \ar\leq\ar
|Z_{t\wedge T_0^x}^x - Z_{t\wedge T_0^x}^y| + \sup_{s\in [T_0^x, T_0^y)} \rho(0,Z_s^y) \cr
 \ar\leq\ar
|x-y|+\sup_{s\in [T_0^x, T_0^y)} Z_s^y.
 \eeqnn
By the right-continuity of the L\'{e}vy process we have a.s.\ $\lim_{t\downarrow T_0^x} Z^y_t =Z^y_{T_0^x}=y-x$. Since the L\'{e}vy process has no  negative jump, we have a.s.\ $\lim_{t\uparrow T_0^y} Z^y_t =0$. Then a.s.\
 $$
\lim_{y\to x} d_\infty(Y^y,Y^x) = \lim_{x\to y} d_\infty(Y^y,Y^x)= 0.
 $$
By Proposition~\ref{t1.0} the process $X^x:= J_\theta(Y^x)$ is a solution to (\ref{e0.1}). A modification of the proof of Proposition 5 in Caballero et al.\ (2009) shows that the transformation $J_\theta$ is continuous on $(D,d_{\infty})$. Then we have a.s.\
 $$
\lim_{y\to x} d_\infty(X^y,X^x) = \lim_{x\to y} d_\infty(X^y,X^x)= 0.
 $$
That proves the desired result. \qed

\medskip

\noindent\emph{Proof of Theorem}~\ref{t1.4}.~ Let $X = J_\theta(Y)$ be constructed as in Proposition~\ref{t1.0}. Let $T_y=\inf\{t\geq 0: Y_t= y\}$. By Corollary~3.13 in Kyprianou (2006, p.81) we have $\mbf{P}(T_y< \infty|Y_0=x)= e^{-q(x-y)}$ for $0\le y\le x$. Since $Y= (Y_t)_{t\ge 0}$ has no negative jumps, we have $Y_t> y$ for $0\le t<T_y$. By Proposition~\ref{t1.0}, for $0< y\le x$ we have $\tau_y= \alpha(T_y)< \infty$ if and only if $T_y< \infty$. Then (\ref{t1.4sss}) holds. Using Proposition~\ref{t1.0} again we see that $\lim_{t\to \infty} X_t= 0$ if and only if $T_0< \infty$. By Proposition 37.10 in Sato (1999, p.255), on the event $\{T_0=\infty\}$ we have a.s.\ $\lim_{t\to \infty} Y_t = \infty$ and hence a.s.\ $\lim_{t\to \infty} X_t = \infty$. Then (\ref{2.10}) holds. \qed

\section{Mean extinction and explosion times}

\setcounter{equation}{0}

In this section we prove the results on the hitting times of the continuous-state polynomial branching process. We shall see that the relations established in Theorem~\ref{t4.1} play important roles in the proofs. Recall that $e_\lambda(x)= e^{-\lambda x}$ for $\lambda\in (0,\infty)$ and $x\in [0,\infty]$.

\bproposition\label{l1x.1xx} \rm(1)~$t\mapsto P_te_\lambda(x)$ is decreasing if $0<\lambda\le q$; \rm(2)~$t\mapsto P_te_\lambda(x)$ is increasing if $q\le \lambda< \infty$; \rm(3)~$\lim_{t\to \infty}P_te_\lambda(x)= e_q(x)$ for all $0< \lambda< \infty$; \rm(4)~$e_q$ is an invariant function of $(P_t)_{t\ge 0}$. \eproposition

\proof It is easy to see that $e_\lambda\in \mcr{D}(L)$ and $Le_\lambda(x) = x^\theta \psi(\lambda) e_\lambda (x)$. By (\ref{d1}) we have
 \beqlb\label{2.7}
\frac{d}{dt}P_te_\lambda(x)= \psi(\lambda) \int_{[0,\infty)} y^\theta e^{-\lambda y} P_t(x,dy).
 \eeqlb
By convexity of $\psi$ we see (1) and (2) hold. By Theorem~\ref{t1.4}\,(2) we get (3), from which (4) follows. \qed

\medskip

\noindent\emph{Proof of Theorem}~\ref{t4.1}.~ Taking the Laplace transform in both sides of (\ref{2.7}) and using integration by parts we get
 \beqnn
\psi(\lambda)\int_{[0,\infty)} y^\theta e^{-\lambda y} U^\eta(x,dy)
 \ar=\ar
\int^\infty_0 e^{-\eta t} \frac{d}{dt}P_te_\lambda(x)dt \cr
 \ar=\ar
e^{-\eta t}P_te_\lambda(x)\big|^{t=\infty}_{t=0} + \eta \int_0^\infty e^{-\eta t}P_t e_\lambda(x)dt \cr
 \ar=\ar
-e^{-\lambda x} + \eta \int_0^\infty e^{-\eta t}P_te_\lambda(x)dt.
 \eeqnn
Then we get~(\ref{2.11}). By~(\ref{de0.2}), we see that
 \beqlb\label{re1}
l_x(\eta,z)
 =
\int_{[0,\infty)} y^\theta e^{-zy} U^\eta(x,dy)
 =
\int^\infty_0 e^{-\eta t}dt\int_{[0,\infty)} y^\theta e^{-zy} P_t(x,dy).
 \eeqlb
Multiplying ({\ref{re1}}) by $(z-\lambda)^{\theta-1}$ and integrating both sides, we have
 \beqnn
\int^\infty_\lambda l_x(\eta,z)(z-\lambda)^{\theta-1} dz
 \ar=\ar
\int^\infty_\lambda(z-\lambda)^{\theta-1} dz\int^\infty_0 e^{-\eta t} dt\int_{[0,\infty)} y^\theta e^{-zy} P_t(x,dy) \cr
 \ar=\ar
\int^\infty_0z^{\theta-1} dz\int^\infty_0 e^{-\eta t}dt\int_{[0,\infty)} y^\theta e^{-(z+\lambda)y} P_t(x,dy) \cr
 \ar=\ar
\int^\infty_0z^{\theta-1} dz\int^\infty_0 e^{-\eta t}dt\int_{(0,\infty)} y^\theta e^{-(z+\lambda)y} P_t(x,dy) \cr
\ar=\ar
\int^\infty_0 e^{-\eta t}dt\int_{(0,\infty)} e^{-\lambda y} P_t(x,dy)\int^\infty_0 {y^\theta z^{\theta-1}} e^{-yz}dz \cr
 \ar=\ar
\Gamma(\theta) \int^\infty_0 e^{-\eta t}dt\int_{(0,\infty)}e^{-\lambda y} P_t(x,dy).
 \eeqnn
That gives (\ref{1.11xx}). \qed

\blemma\label{l1.1} For any $\lambda\ge 0$ and $x\in [0,\infty)$ we have
 \beqlb\label{4.3}
e^{-qx}-e^{-\lambda x} = \psi(\lambda) \int^\infty_0 dt\int_{[0,\infty)} y^\theta e^{-\lambda y} P_t(x,dy)
 \eeqlb
and
 \beqlb\label{1.11}
\int^\infty_\lambda h_x(z)(z-\lambda)^{\theta-1} dz
=
\Gamma(\theta)\int^\infty_0 dt\int_{(0,\infty)}e^{-\lambda y} P_t(x,dy).
 \eeqlb

\elemma

\proof By Proposition~ \ref{t1.4sss} and Proposition~\ref{l1x.1xx}, as $t\uparrow \infty$ we have $e^{-\lambda x}\le P_te_\lambda(x)\uparrow e^{-qx}$ for $q< \lambda< \infty$ and $e^{-\lambda x}\ge P_te_\lambda(x)\downarrow e^{-qx}$ for $0< \lambda< q$. Since
 \beqnn
\eta U^\eta e_\lambda(x)
 =
\int_0^\infty \eta e^{-\eta t}P_te_\lambda(x) dt
 =
\int_0^\infty e^{-t}P_{t/\eta}e_\lambda(x) dt,
 \eeqnn
we see that $e^{-\lambda x}\le \eta U^\eta e_\lambda(x)\uparrow e^{-qx}$ for $q< \lambda< \infty$ and $e^{-\lambda x}\ge \eta U^\eta e_\lambda(x)\downarrow e^{-qx}$ for $0< \lambda< q$ as $\eta\downarrow 0$. Then we use monotone convergence to get (\ref{4.3}) and (\ref{1.11}) by letting $\eta\to 0$ in (\ref{2.11}) and (\ref{1.11xx}), respectively. \qed

\medskip

\noindent\emph{Proof of Theorem}~\ref{t1.7}.~ Observe that, for any $\eta\ge 0$,
 \beqlb\label{4.5vv}
\int_{(0,\infty)}e^{-\eta z} P_t(x,dz)
 \ar=\ar
\frac{1}{\Gamma(\theta)}\int_{(0,\infty)}e^{-\eta z} P_t(x,dz)\int^\infty_0 y^{\theta-1} e^{-y}dy \cr
 \ar=\ar
\frac{1}{\Gamma(\theta)}\int_{(0,\infty)} z^\theta e^{-\eta z} P_t(x,dz)\int^\infty_0 \lambda^{\theta-1} e^{-\lambda z}d\lambda \cr
 \ar=\ar
\frac{1}{\Gamma(\theta)}\int^\infty_0\lambda^{\theta-1}d\lambda \int_{(0,\infty)} z^\theta e^{-(\lambda+\eta)z} P_t(x,dz).
 \eeqlb
By Theorem~\ref{t1.4}\,(2) and Proposition~\ref{l1x.1xx}\,(4), the function $x\mapsto e^{-qx}= \mbf{P}_x(X_\infty= 0)= \mbf{P}_x(\tau_\infty> t, X_\infty= 0)$ is invariant for the transition semigroup of $X$. By (\ref{4.5vv}),
 \beqnn
\mbf{P}_x(\tau_0> t, X_\infty= 0)
 \ar=\ar
\mbf{P}_x(X_\infty= 0) - \mbf{P}_x(\tau_0\leq t, X_\infty= 0) \\
 \ar=\ar
e^{-qx} - P_t(x,\{0\})
 =
\int_{(0,\infty)}e^{-qz} P_t(x,dz) \cr
 \ar=\ar
\frac{1}{\Gamma(\theta)}\int^\infty_0\lambda^{\theta-1}d\lambda \int_{(0,\infty)} z^\theta e^{-(\lambda+q)z} P_t(x,dz)
 \eeqnn
and
 \beqnn
\mbf{P}_x(\tau_\infty> t, X_\infty= \infty)
 \ar=\ar
\mbf{P}_x(\tau_\infty> t) - \mbf{P}_x(\tau_\infty> t, X_\infty= 0) \\
 \ar=\ar
P_t(x,[0,\infty)) - e^{-qx}
 =
\int_{(0,\infty)}(1-e^{-qz}) P_t(x,dz) \cr
 \ar=\ar
\frac{1}{\Gamma(\theta)}\int^\infty_0\lambda^{\theta-1}d\lambda \int_{(0,\infty)} z^\theta e^{-\lambda z}(1-e^{-qz}) P_t(x,dz).
 \eeqnn
By Lemma~\ref{l1.1} we have
 \beqnn
\mbf{E}_x(\tau_0: X_\infty= 0)
 \ar=\ar
\int^\infty_0 \mbf{P}_x(\tau_0> t, X_\infty= 0)dt \cr
 \ar=\ar
\frac{1}{\Gamma(\theta)}\int^\infty_0\lambda^{\theta-1}d\lambda \int^\infty_0 dt \int_{(0,\infty)} z^\theta e^{-(\lambda+q)z} P_t(x,dz) \cr
 \ar=\ar
\frac{1}{\Gamma(\theta)}\int^\infty_0 h_x(\lambda+q)\lambda^{\theta-1} d\lambda
 \eeqnn
and
 \beqnn
\mbf{E}_x(\tau_\infty: X_\infty= \infty)
 \ar=\ar
\int_0^\infty \mbf{P}_x(\tau_\infty> t, X_\infty= \infty) dt \cr
 \ar=\ar
\frac{1}{\Gamma(\theta)}\int^\infty_0\lambda^{\theta-1}d\lambda \int^\infty_0 dt \int_{(0,\infty)} z^\theta e^{-\lambda z}(1-e^{-qz}) P_t(x,dz) \cr
 \ar=\ar
\frac{1}{\Gamma(\theta)}\int^\infty_0 [h_x(\lambda) - h_x(\lambda+q)] \lambda^{\theta-1}d\lambda,
 \eeqnn
where $h_x$ is defined by~(\ref{de0.3}). Then (\ref{2.3}) and (\ref{2.3'}) are proved. By summing up those two expressions we get (\ref{2.4}). \qed

\bproposition\label{t4.1ss} For any $y\in [0,\infty)$, let $(P_t^{(y)})_{t\ge 0}$ denote the transition semigroup on $[y,\infty]$ of the stopped process $X^{(y)}= (X_{t\land \tau_y}: t\ge 0)$ with $X_0\ge y$. For any $\eta>0$, $\lambda\ge 0$ and $x\in [y,\infty)$ we have
 \beqlb\label{4.6}
\eta\int_0^\infty e^{-\eta t}P_t^{(y)}e_\lambda(x)dt
 =
e^{-\lambda x} + \psi(\lambda)\int_0^\infty e^{-\eta t}dt\int_{(y,\infty)} z^\theta e^{-\lambda z}P_t^{(y)}(x,dz).
 \eeqlb
\eproposition

\proof For any $x\ge y\in [0,\infty)$ we have $\mbf{P}_x(\tau_y< \tau_0, X_{\tau_y}= y)= 1$. Then (\ref{2.6}) implies
 \beqnn
e^{-\lambda X_{t\land \tau_y}} \ar=\ar e^{-\lambda x} + \psi(\lambda)\int^{t\land \tau_y}_0 X_s^\theta e^{-\lambda X_s}ds + \mbox{martingale} \cr
 \ar=\ar
e^{-\lambda x} + \psi(\lambda)\int^t_0 X_{s\land \tau_y}^\theta e^{-\lambda X_{s\land \tau_y}} 1_{(y,\infty)}(X_{s\land \tau_y})ds + \mbox{martingale}.
 \eeqnn
Taking the expectation in both sides yields
 \beqnn
P_t^{(y)}e_\lambda(x) = e^{-\lambda x} + \psi(\lambda)\int^t_0 ds \int_{(y,\infty)}z^\theta e^{-\lambda z} P_s^{(y)}(x,dz).
 \eeqnn
Thus we have
 \beqnn
\eta\int_0^\infty e^{-\eta t}P_t^{(y)}e_\lambda(x)dt
 \ar=\ar
e^{-\lambda x} + \eta\psi(\lambda)\int_0^\infty e^{-\eta t}dt\int_0^tds\int_{(y,\infty)} z^\theta e^{-\lambda z}P_s^{(y)}(x,dz) \cr
 \ar=\ar
e^{-\lambda x} + \eta\psi(\lambda)\int_0^\infty ds\int_s^\infty e^{-\eta t}dt\int_{(y,\infty)} z^\theta e^{-\lambda z}P_s^{(y)}(x,dz) \cr
 \ar=\ar
e^{-\lambda x} + \eta\psi(\lambda)\int_0^\infty e^{-\eta s}ds\int_0^\infty e^{-\eta t}dt\int_{(y,\infty)} z^\theta e^{-\lambda z}P_s^{(y)}(x,dz) \cr
 \ar=\ar
e^{-\lambda x} + \psi(\lambda)\int_0^\infty e^{-\eta t}dt\int_{(y,\infty)} z^\theta e^{-\lambda z}P_t^{(y)}(x,dz).
 \eeqnn
That proves (\ref{4.6}). \qed

\medskip

\noindent\emph{Proof of Theorem}~\ref{t1.8}.~
By Theorem~\ref{t1.4}\,(2) we have $\mbf{P}_x(\{X_\infty= 0\}\cup \{X_\infty= \infty\})= 1$. Observe that
\beqnn
\eta\int^\infty_0 e^{-\eta t}P_t^{(y)}e_\lambda(x) dt
 \ar=\ar
\int^\infty_0 e^{-t} dt\int^\infty_y e^{-\lambda z} P^{(y)}_{t/\eta}(x, dz)\cr
 \ar=\ar
\int^\infty_0 e^{-t}\mbf{E}_x\left[\exp\left\{-\lambda X_{(t/\eta)\wedge \tau_y}\right\} 1_{\{\tau_y<t/\eta\}}\right] dt \cr
 \ar\ar
+ \int^\infty_0 e^{-t}\mbf{E}_x\left[\exp\left\{-\lambda X_{(t/\eta)\wedge \tau_y}\right\} 1_{\{\tau_y\ge t/\eta\}}\right] dt.
\eeqnn
 The first term on the right-hand side converges to $e^{-\lambda y}\mbf{P}_x(\tau_y< \infty)$ as $\eta\to 0$. Since the process $X$ started from $x\,(> y)$ can come to $0$ only by crossing $y$, we have a.s.\ $X_\infty = \infty$ on the event $\{\tau_y = \infty\}$, so the second term vanishes as $\eta\to 0$. Then by Theorem~\ref{t1.4}\,(1), as $\eta\to 0$ we have
 $$
\eta\int^\infty_0 e^{-\eta t}P_t^{(y)}e_\lambda(x)dt
\to
e^{-\lambda y-q(x-y)}.
 $$
By taking $\eta\to 0$ in (\ref{4.6}) we obtain
 $$
e^{-qx}(e^{-(\lambda-q)y}-e^{-(\lambda-q)x})
 =
\psi(\lambda) \int^\infty_0 dt\int_{(y,\infty)} z^\theta e^{-\lambda z} P_t^{(y)}(x,dz).
 $$
Thus we have
 \beqnn
e^{-qx}\int^\infty_0\frac{e^{-(\lambda-q)y}-e^{-(\lambda-q)x}}{\psi(\lambda)}\lambda^{\theta-1} d\lambda
 \ar=\ar
\int^\infty_0\lambda^{\theta-1} d\lambda\int^\infty_0 dt\int_{(y,\infty)} z^\theta e^{-\lambda z} P_t^{(y)}(x,dz) \cr
 \ar=\ar
\int^\infty_0\eta^{\theta-1}e^{-\eta} d\eta\int^\infty_0 dt\int_{(y,\infty)} P_t^{(y)}(x,dz) \cr
 \ar=\ar
\Gamma(\theta)\int^\infty_0 P_t^{(y)}(x,(y,\infty))dt \cr
 \ar=\ar
\Gamma(\theta)\int^\infty_0 \mbf{P}_x(\tau_\infty\land \tau_y> t)dt.
 \eeqnn
That implies (\ref{0.1''}) by a formula for the expectation. \qed

\medskip

\noindent\emph{Proof of Corollary}~\ref{t1.6'}.~ By Theorem~{\ref{t1.4}\,(1) we see $\mbf{P}_x(\tau_y< \infty)= 1$ and hence $\mbf{P}_x(\tau_y< \tau_\infty)= 1$. Then~(\ref{2.9}) follows from~(\ref{0.1''}). \qed

\section{Extinction and explosion probabilities}

\setcounter{equation}{0}

In this section we give the proofs of the results on the extinction and explosion probabilities of the continuous-state polynomial branching process.

\blemma\label{l1.4a} Let $\varepsilon>q$. Then for any $x\in (0,\infty)$, we have $\mbf{P}_x(\tau_0< \infty)> 0$ if and only if
 (\ref{1.1'})
holds.
 \elemma

\proof (1) Suppose that $\mbf{P}_x(\tau_0< \infty)> 0$ for some $x\in (0,\infty)$. Then for sufficiently large $t\ge 0$ we have $P_t(x,\{0\})= \mbf{P}_x(\tau_0\le t)> 0$. It follows that
 $$
\rho_x(\eta) := \downarrow\!\!\lim_{\lambda\uparrow\infty} U^\eta e_\lambda(x)= \int^\infty_0 e^{-\eta t} P_t(x, \{0\})dt> 0, \qquad \eta> 0.
 $$
Then there exists $\varepsilon_0= \varepsilon_0(x,\eta)> q$ such that $e^{-\varepsilon_0x}\le \eta\rho_x(\eta)/2$, and so
 $$
\eta U^\eta e_\lambda(x) - e^{-\lambda x}
 \ge
\eta U^\eta e_\lambda(x) - \frac{1}{2}\eta\rho_x(\eta)
 \ge
\frac{1}{2}\eta\rho_x(\eta), \qquad \lambda> \varepsilon_0.
 $$
For $\lambda\ge \varepsilon_0+1$ we have $\lambda^{\theta-1}\le (\lambda-\varepsilon_0)^{\theta-1}$ if $0<\theta\le 1$, and $\lambda^{\theta-1}\le (\varepsilon_0 + 1)^{\theta-1}(\lambda-\varepsilon_0)^{\theta-1}$ if $\theta> 1$. Therefore we can find a constant $C=C(\varepsilon_0,\theta)> 0$ such that
 \beqnn
\int^\infty_{\varepsilon_0+1} \frac{\lambda^{\theta-1}}{\psi(\lambda)} d\lambda
 \ar\le\ar
\frac{C}{2}\int^\infty_{\varepsilon_0} \frac{(\lambda-\varepsilon_0)^{\theta-1}}{\psi(\lambda)} d\lambda
 \le
\frac{C}{\eta\rho_x(\eta)}\int^\infty_{\varepsilon_0} \frac{\eta U^\eta e_\lambda(x)-e^{-\lambda x}}{\psi(\lambda)} (\lambda-\varepsilon_0)^{\theta-1} d\lambda.
 \eeqnn
By (\ref{1.11xx}) we see the right hand side is finite, and hence~(\ref{1.1'}) holds for any $\varepsilon> q$.}

(2) Suppose that $\mbf{P}_x(\tau_0< \infty)= 0$ for some $x\in (0,\infty)$. Then $P_t(x,\{0\})= 0$ for every $t\ge 0$. By Theorem~\ref{t1.4}\,(2), for $\varepsilon>q$ we have
 $$
\lim_{t\to \infty}\int_{(0,\infty)}e^{-\varepsilon z} P_t(x,dz)
 =
\lim_{t\to \infty}\int_{[0,\infty)}e^{-\varepsilon z} P_t(x,dz)= e^{-qx}>0.
 $$
By (\ref{1.11}) it follows that
 $$
\int^\infty_\varepsilon \frac{\lambda^{\theta-1}}{\psi(\lambda)}d\lambda
 \geq
\int^\infty_\varepsilon h_x(\lambda)(\lambda-\varepsilon)^{\theta-1}d\lambda
 =
\Gamma(\theta)\int^\infty_0 dt\int_{(0,\infty)}e^{-\varepsilon z} P_t(x,dz)= \infty.
 $$\qed

\blemma\label{l1.4b}
Suppose that $\mbf{P}_x(\tau_0< \infty)> 0$ for some $x\in (0,\infty)$. Then we have $\mbf{P}_x(\tau_0< \infty)= \mbf{P}_x(X_\infty= 0)= e^{-qx}.$
\elemma

\proof Suppose that $\mbf{P}_x(\tau_0< \infty)> 0$ for some $x\in (0,\infty)$. By Theorem~\ref{t1.4}\,(2) we only need to prove $\mbf{P}_x(\tau_0= \infty, X_\infty= 0)= 0$. In this case, we have $\mbf{P}_x(\tau_0= \infty)< 1$, and hence $\alpha:= \mbf{P}_x (\tau_0> v)= \mbf{P}_x (X_v>0)< 1$ for some $v> 0$. By Theorem~\ref{t1.1} we have $\mbf{P}_y(X_v>0)< \alpha$ for $y\le x$. Let $\sigma_0=0$ and $\sigma_n=\inf\{t> \sigma_{n-1}+v: X_t\leq x\}$ for $n\ge 1$. It is easy to see that $X_{\sigma_n}\le x$. By the strong Markov property, for any $n\ge 1$ we have
 \beqnn
\mbf{P}_x\big(\tau_0= \infty, X_\infty= 0\big)
 \ar\le\ar
\mbf{P}_x\bigg(\bigcap_{k=1}^n \big\{\sigma_k< \infty, X_{\sigma_k+v}> 0\big\}\bigg) \cr
 \ar\le\ar
\mbf{E}_x\bigg[\prod_{k=1}^{n-1} 1_{\{\sigma_k< \infty, X_{\sigma_k+v}> 0\}}1_{\{\sigma_n< \infty\}} \mbf{P}_x\big(X_{\sigma_n+v}> 0|\mcr{F}_{\sigma_n}\big)\bigg] \cr
 \ar\le\ar
\mbf{E}_x\bigg[\prod_{k=1}^{n-1} 1_{\{\sigma_k< \infty, X_{\sigma_k+v}> 0\}}1_{\{\sigma_n< \infty\}} \mbf{P}_{X_{\sigma_n}}(X_v> 0)\bigg] \cr
 \ar\le\ar
\alpha\mbf{E}_x\bigg[\prod_{k=1}^{n-1} 1_{\{\sigma_k< \infty, X_{\sigma_k+v}> 0\}}1_{\{\sigma_n< \infty\}}\bigg] \cr
 \ar\le\ar
\alpha\mbf{P}_x\bigg(\bigcap_{k=1}^{n-1} \{\sigma_k< \infty, X_{\sigma_k+v}> 0\}\bigg)
 \le
\dots\le \alpha^n.
 \eeqnn
Then the left-hand side vanishes. \qed

\medskip

\noindent\emph{Proof of Theorem}~\ref{t1.4'}.~ From Lemmas~\ref{l1.4a} and \ref{l1.4b}, we only need to check that (\ref{1.1'}) does not hold in the case $\theta\ge 2$. By the Taylor expansion, we see $e^{-\lambda u} - 1 + \lambda u\le \lambda^2u^2/2$. In view of (\ref{0.1}) we have
 \beqnn
\psi(\lambda) \le |b|\lambda + c\lambda^2 + \frac{1}{2} \lambda^2 \int_{(0,1]} u^2 m(du) - \lambda\int_{(1,\infty]} (1-e^{-\lambda u}) m(du).
 \eeqnn
Then there is a constant $C>0$ so that $\psi(\lambda)\le C\lambda^2$ for $\lambda\geq \varepsilon$. If $\theta\geq 2$, then
 $$
\int^{\infty}_\varepsilon \frac{\lambda^{\theta-1}}{\psi(\lambda)} d\lambda
 \ge
\frac{1}{C}\int^{\infty}_\varepsilon \lambda^{\theta-3} d\lambda= \infty,
 $$
so (\ref{1.1'}) does not hold. \qed

\medskip

\noindent\emph{Proof of Corollary}~\ref{r1.4'}.~ If $c>0$, we can take $\varepsilon>0$ so that $\psi(\lambda)\ge c\lambda^2/2$ for $\lambda\ge \varepsilon$. When $0<\theta<2$, we have
 $$
\int^{\infty}_\varepsilon \frac{\lambda^{\theta-1}}{\psi(\lambda)} d\lambda
 \le
\frac{2}{c}\int^{\infty}_\varepsilon \lambda^{\theta-3} d\lambda<\infty,
 $$
so the process hits $0$ by Lemma~\ref{l1.4a}. \qed

\bproposition\label{l1.3a}
Let $\varepsilon> q$. For any $x\in (0,\infty)$ we have $\mbf{P}_x(\tau_\infty< \infty)= 0$ if and only if $a= -\psi(0)= 0$ and one of the following two conditions is satisfied: {\rm(i)} $\psi'(0)\geq 0$; {\rm(ii)} $\psi'(0)< 0$ and (\ref{1.1}) holds. \eproposition

\proof (1) In the case $a= -\psi(0)> 0$, we can let $\lambda\to 0$ in~(\ref{2.11}) to see
 \beqnn
\eta\int^\infty_0 e^{-\eta t}P_t(x,[0,\infty))dt
 =
\eta U^\eta(x,[0,\infty))
 =
1 + \psi(0)\int_{[0,\infty)} z^\theta U^\eta(x,dz)< 1.
 \eeqnn
Then for some $t>0$ we have $P_t(x,[0,\infty))< 1$ and so $\mbf{P}_x (\tau_\infty\le t)= P_t(x,\{\infty\})> 0$.

(2) Suppose that $a= -\psi(0)= 0$ and $\psi'(0)\geq 0$. By the convexity of $\psi$ we have $\psi(\lambda)>0$ for each $\lambda> 0$. Then~(\ref{2.11}) implies
 $$
\eta\int_{[0,\infty)}e^{-\lambda z}U^\eta(x,dz)
 =
\eta\int^\infty_0 e^{-\eta t}dt\int_{[0,\infty)}e^{-\lambda z}P_t(x,dz)> e^{-\lambda x}.
 $$
By letting $\lambda\to 0$ on the both sides we see $\eta U^\eta(x,{[0,\infty)})= 1$. Then $\mbf{P}_x (\tau_\infty> t)= P_t(x,[0,\infty))= 1$ for every $t>0$. That implies $\mbf{P}_x(\tau_\infty= \infty)= 1$.

(3) Consider the case with $a= -\psi(0)= 0$ and $\psi'(0)< 0$. (a) Suppose that~(\ref{1.1}) holds but $\mbf{P}_x(\tau_\infty< \infty)> 0$. Then $P_t(x,[0,\infty))= \mbf{P}_x(\tau_\infty> t)< 1$ for sufficiently large $t\ge 0$. For any $\eta>0$ we have
 \beqlb\label{1.3}
\kappa(x):= 1 - \eta\int^\infty_0 e^{-\eta t}P_t(x, [0,\infty))dt> 0.
 \eeqlb
By continuity there exists an $\varepsilon\in (0,q)$ such that $\psi(\lambda)<0$ and
 \beqnn
e^{-\lambda x}-\eta U^\eta e_\lambda(x)\geq \frac{1}{2}\kappa(x)> 0, \qquad 0<\lambda\le \varepsilon.
 \eeqnn
By~(\ref{1.1}) we have
 $$
\int^\varepsilon_0 l_x(\eta,\lambda)\lambda^{\theta-1}d\lambda
 =
\int^\varepsilon_0 \frac{\eta U^\eta e_\lambda(x)-e^{-\lambda x}}{\psi(\lambda)} \lambda^{\theta-1}d\lambda
 \geq
\frac{\kappa(x)}{2} \int^\varepsilon_0 \frac{\lambda^{\theta-1}}{-\psi(\lambda)}d\lambda=\infty.
 $$
Then~(\ref{1.11xx}) implies
 \beqnn
\int^\infty_0 e^{-\eta t} P_t(x, (0,\infty))dt
 =
\frac{1}{\Gamma(\theta)}\int^\infty_0 l_x(\eta,\lambda)\lambda^{\theta-1} d\lambda= \infty,
 \eeqnn
which is in contradiction to~(\ref{1.3}). (b) Conversely, suppose that~(\ref{1.1}) does not hold. Then we have
 \beqnn
 \int^\varepsilon_0 \frac{\lambda^{\theta-1}}{-\psi(\lambda)}d\lambda< \infty.
 \eeqnn
Using the convexity of $\psi$ we know $\psi'(q)> 0$, and so
 $$
\lim_{\lambda\to q} \frac{e^{-(\lambda-q)y}-e^{-(\lambda-q)x}}{\psi(\lambda)}
 =
\frac{(x-y)}{\psi'(q)}.
 $$
Since $\lim_{\lambda\to\infty}\psi(\lambda)=\infty$, by Theorem~\ref{t1.8} we see
 \beqnn
\mbf{E}_x(\tau_\infty: \tau_y= \infty)
 \le
\mbf{E}_x(\tau_\infty\wedge \tau_y)
 =
\frac{\,e^{-qx}}{\Gamma(\theta)}\int_0^\infty\frac{e^{-(\lambda-q)y}-e^{-(\lambda-q)x}} {\psi(\lambda)} \lambda^{\theta-1}d\lambda< \infty.
 \eeqnn
It follows that $\mbf{P}_x(\tau_\infty< \infty)\ge \mbf{P}_x(\tau_y= \infty)= 1-e^{-q(x-y)}> 0$.\qed

\blemma\label{l1.3b}
Suppose that $\mbf{P}_x(\tau_\infty< \infty)> 0$ for some $x\in(0,\infty)$. Then we have $\mbf{P}_x (\tau_\infty< \infty) = \mbf{P}_x(X_\infty= \infty) = 1-e^{-qx}.$
\elemma
\proof We only need to prove $\mbf{P}_x(\tau_\infty= \infty, X_\infty= \infty)= 0$ for each $x\in (0,\infty)$. Fix $x\in (0,\infty)$ and choose sufficiently large $v>0$ so that $\alpha:= \mbf{P}_x (\tau_\infty > v)= \mbf{P}_x (X_v< \infty)< 1$. By Theorem~\ref{t1.1} we have $\mbf{P}_y(X_v< \infty)\le \alpha$ for $y\ge x\in [0,\infty)$. Let $\sigma_1=0$ and $\sigma_n=\inf\{t> \sigma_{n-1}+v: X_t\geq x\}$ for $n\ge 1$. As in the proof of Theorem~\ref{t1.4'} one sees
 \beqnn
\mbf{P}_x\big(\tau_\infty= \infty, X_\infty= \infty\big)
 \le
\mbf{P}_x\bigg(\bigcap_{k=1}^n \big\{\sigma_k< \infty, X_{\sigma_k+v}< \infty\big\}\bigg)
 \le
\alpha^n.
 \eeqnn
for every $n\ge 1$. Then we must have $\mbf{P}_x(\tau_\infty= \infty, X_\infty= \infty)= 0$. \qed

\medskip

\noindent\emph{Proof of Theorem}~\ref{t1.3}.~ (1) Suppose that $\theta>1$. By  Proposition~\ref{l1.3a}, we have $\mbf{P}_x(\tau_\infty< \infty)> 0$ if $a= -\psi(0)> 0$, and $\mbf{P}_x(\tau_\infty< \infty)= 0$ if $\psi(0)=0$ and $\psi'(0)\ge 0$. Since $\theta> 1$, when $a= -\psi(0)= 0$ and $\psi'(0)< 0$, we have
 $$
\int_0^\varepsilon\frac{\lambda^{\theta-1}}{-\psi(\lambda)}d\lambda< \infty.
 $$
Then $\mbf{P}_x(\tau_\infty< \infty)> 0$ by Proposition~\ref{l1.3a}.

(2) Suppose that $0<\theta\le 1$. It suffices to consider the case with $a= -\psi(0)= 0$ and $\psi'(0)>-\infty$. In this case, since $0<\theta\leq 1$ and $\psi(\lambda)= \psi'(0)\lambda + o(\lambda)$ as $\lambda\to 0$, we have
 $$
\int_0^\varepsilon\frac{\lambda^{\theta-1}}{-\psi(\lambda)}d\lambda= \infty.
 $$
Then $\mbf{P}_x(\tau_\infty< \infty)= 0$ by Proposition~\ref{l1.3a}. Finally, by using Lemma~\ref{l1.3b} we complete the proof. \qed

\section{The process coming down from infinity}

\setcounter{equation}{0}

In this section, we give a construction of the continuous-state polynomial branching process coming down from $\infty$.

\medskip

\noindent\emph{Proof of Theorem}~\ref{t3.2}.~By the right continuity of $(X^x_t)_{t\ge 0}$ we have a.s.\ $\tau^x_y>0$ for $x> y\in [0,\infty)$, yielding a.s.\ $\tau^\infty_y> 0$ for $y\in [0,\infty)$. Then (iii) $\Rightarrow$ (ii) $\Rightarrow$ (i). By~(\ref{1.5}) we see (iii) $\Leftrightarrow$ (iv). To show that (i) $\Rightarrow$ (iii), suppose that $\mbf{P}(\tau^{\infty}_y< \infty)> 0$ for some $y>0$. Then there exists $t>0$ such that $\alpha:= \mbf{P}(\tau^{\infty}_y> t)< 1$. By Theorem~\ref{t1.1} we see for each $x\ge y\in [0,\infty)$ we have $\mbf{P}(\tau^x_y> t)\leq\alpha< 1$.
By the Markov property, for $n\ge 1$,
 \beqnn
\mbf{P}(\tau^\infty_y> nt)
\ar=\ar
\lim_{x\to\infty}\mbf{P}(\tau^x_y>nt) \\
\ar=\ar
\lim_{x\to\infty}\mbf{E}(1_{\{\tau^x_y> t\}}1_{\{\tau^x_y\circ\theta_t> (n-1)t\}})\\
 \ar=\ar
\lim_{x\to\infty}\mbf{E}[1_{\{\tau^x_y> t\}}\mbf{E}(1_{\{\tau^x_y \circ\theta_t> (n-1)t\}}|\mcr{F}_t)]\\
 \ar=\ar
\lim_{x\to\infty}\mbf{E}[1_{\{\tau^x_y> t\}}\mbf{E}(1_{\{\tau^z_y> (n-1)t\}})|_{z=X^x_t}]\\
 \ar\le\ar
\mbf{E}[1_{\{\tau^\infty_y> t\}}\mbf{E}(1_{\{\tau^\infty_y> (n-1)t\}})] \\
 \ar=\ar
\alpha\mbf{P}(\tau^\infty_y> (n-1)t)].
 \eeqnn
Then $\mbf{P}(\tau^\infty_y> nt)\leq \alpha^n$ by induction. That implies (iii). \qed

\blemma\label{l5.1}
For any $0<x<\infty$ we have a.s.\ $\lim_{y\uparrow x}\tau_0^y=\tau^x_0$.
\elemma

\proof By Theorem~\ref{t1.1} the mapping $y\mapsto \tau^y_0$ is increasing. Then for $x>y>0$ and $\lambda>0$ we have
 $$
\mbf{E}[e^{-\lambda \tau^y_0}]\ge \mbf{E}[e^{-\lambda \tau^x_0}].
 $$
By Fubini theorem, for $\lambda, y>0$,
 \beqnn
\mbf{E}[e^{-\lambda \tau^y_0}]
 =
\mbf{E}\bigg[1-\lambda\int^\infty_0 e^{-\lambda t}1_{\{t<\tau^y_0\}}dt\bigg]
 =
1-\lambda\int^\infty_0 e^{-\lambda t}\mbf{P}(X^y_t>0)dt.
\eeqnn
By Theorem~\ref{t1.4xx}, we infer $X^y_t\to X^x_t$ in distribution as $y\to x$. Then, by the above equality and Theorem 3.1~(e) in Ethier and Kurtz (1986, page 108),
\beqnn
\lim_{y\uparrow x}\mbf{E}[e^{-\lambda \tau^y_0}]\le\mbf{E}[e^{-\lambda \tau^x_0}].
\eeqnn
From the above two inequalities it follows that $\tau^y_0\to \tau^x_0$ in distribution as $y\uparrow x$. By the monotonicity of $y\mapsto \tau^y_0$ we see a.s.\ $\lim_{y\uparrow x}\tau_0^y= \tau^x_0$. \qed

\bproposition\label{l1.2}
Suppose that $a= -\psi(0)= 0$ and $\psi'(0)\geq 0$. Then for each $x\in (0,\infty)$ we have a.s.\ $\lim_{y\uparrow x} \sup_{t\ge 0}|X^x_t- X^y_t|= 0$.
\eproposition

\proof
Let $\beta$ be defined as in (\ref{0.2'}). Under the assumptions, we have $\beta>-\infty$ and $a= m(\{\infty\})= 0$, so we can rewritten (\ref{e0.1}) as
\beqnn
X_t\ar=\ar X_0 -\beta \int^{t}_0 X_{s-}^\theta ds+\sqrt{2c}\int^{t}_0 X_{s-}^{\theta/2}dB_s + \int_0^t\int_{(0,\infty)} \int^{X_{s-}^\theta}_0 z \tilde{M}(ds,dz,du).
\eeqnn
For each $n\ge 1$ define the function $r_n$ as in (\ref{0.4}). For each $x>0$ let $\{\xi^x_n(t): t\ge 0\}$ be the unique solution to the following equation
\beqnn
\xi_n(t) \ar=\ar x + \int^t_0 \sqrt{2cr_n(\xi_n(s-))} dB(s) - \int^t_0 \beta r_n(\xi_n(s-))ds \cr
 \ar\ar\qquad
+ \int^t_0 \int_{(0,\infty)}\int^{r_n(\xi_n(s-))}_0 (z\land n) \tilde{M}(ds,dz,du).
 \eeqnn
For $0< y< x$ define $\zeta^{x,y}_n = \inf\{t\geq0: X^x_t\ge n$ or $X^y_t\le 1/n\}.$
Since $1/n \leq X^y_t= \xi_n^y(t)\leq X^x_t= \xi_n^x(t) \leq n$ for $0\leq t< \zeta^{x,y}_n$, the trajectory $t\mapsto \xi_n^x(t)$ and $t\mapsto \xi_n^y(t)$ have no jumps larger than $n$ on the time interval $[0,\zeta^{x,y}_n).$
Then we have
\beqnn
\xi_n^x(t\wedge\zeta^{x,y}_n)-\xi_n^y(t\wedge\zeta^{x,y}_n)
 \ar=\ar
x-y -\beta \int^{t\wedge \zeta^{x,y}_n}_0 [\xi_n^x(s-)^\theta-\xi_n^y(s-)^\theta] ds\cr
 \ar\ar\quad
+\,\sqrt{2c}\int^{t\wedge \zeta^{x,y}_n}_0 [\xi_n^x(s-)^{\theta/2}-\xi_n^y(s-)^{\theta/2}] dB_s\cr
 \ar\ar\quad
+\int_0^{t\wedge \zeta^{x,y}_n}\int_{(0,\infty)} \int^{\xi_n^x(s-)^\theta}_{\xi_n^y(s-)^\theta} (z\land n) \tilde{M}(ds,dz,du).
 \eeqnn
By applying Doob's inequality to the martingale terms and applying H\"{o}lder's inequality to the drift term in above we have
\beqnn
\ar\ar\mbf{E}\bigg[\sup_{0\leq s\leq t\wedge\zeta^{x,y}_n} |\xi_n^x(s-)- \xi_n^y(s-)|^2\bigg] \cr
 \ar\ar\qquad
\leq 4|x-y|^2+ 4\beta^2\mbf{E}\bigg[\bigg(\int^t_0 |\xi_n^x(s-)^\theta-\xi_n^y(s-)^\theta|ds\bigg)^2 \bigg]\cr
 \ar\ar\qquad\qquad
+\, 32c\mbf{E}\bigg[\int^t_0 |\xi_n^x(s-)^{\theta/2}-\xi_n^y(s-)^{\theta/2}|^2ds\bigg]\cr
 \ar\ar\qquad\qquad
+\, 16 \mbf{E}\bigg[\int^t_0 |\xi_n^x(s-)^{\theta}-\xi_n^y(s-)^{\theta}|^2ds\int_{(0,\infty)}z^2\wedge n^2m(dz)\bigg]\cr
 \ar\ar\qquad
\leq 4|x-y|^2+ 4\beta^2t\int^t_0 \mbf{E}[|\xi_n^x(s-)^\theta-\xi_n^y(s-)^\theta|^2]ds\cr
 \ar\ar\qquad\qquad
+\, 32c \int^t_0 \mbf{E}[|\xi_n^x(s-)^{\theta/2}-\xi_n^y(s-)^{\theta/2}|^2]ds\cr
 \ar\ar\qquad\qquad
+\, 16\int^t_0 \mbf{E}[|\xi_n^x(s-)^{\theta}-\xi_n^y(s-)^{\theta}|^2]ds\int_{(0,\infty)}z^2\wedge n^2m(dz).
\eeqnn
Obviously $x\mapsto x^\theta$ and $x\mapsto x^{\theta/2}$ are Lipschitz functions on $[1/n,n]$. Then for $t<k$ there exists a constant $C_{n,k}$ such that
\beqnn
\mbf{E}\bigg[\sup_{0\leq s\leq t\wedge\zeta^{x,y}_n} |\xi_n^x(s)- \xi_n^y(s)|^2\bigg]
 \ar\leq\ar
4|x-y|^2+ C_{n,k}\int^t_0 \mbf{E}(|\xi_n^x(s-)- \xi_n^y(s-)|^2) ds\cr
 \ar\leq\ar
4|x-y|^2+ C_{n,k}\int^t_0 \mbf{E}\bigg[\sup_{0\leq u\leq s\wedge\zeta^{x,y}_n}|\xi_n^x(s)- \xi_n^y(s)|^2\bigg] ds.
\eeqnn
Then by Gronwall's inequality for $t\leq k$ we see
\beqnn
\mbf{E}\bigg[\sup_{0\leq s\leq t\wedge\zeta^{x,y}_n} |\xi_n^x(s)- \xi_n^y(s)|^2\bigg]
\leq
4|x-y|^2\exp[C_{n,k}t].
\eeqnn
It follows that
$$
\lim_{y\uparrow x} \mbf{E} \bigg[\sup_{0\leq s\leq k\land\zeta^{x,y}_n} |\xi_n^x(s)- \xi_n^y(s)|^2\bigg]= 0,
$$
which yields
$$
\lim_{y\uparrow x} \sup_{0\leq s\leq k\land\zeta^{x,y}_n} |\xi_n^x(s)- \xi_n^y(s)|= 0 \qquad\mbox{a.s.}
$$
Since $X^x$ has no negative jump, for any $\varepsilon>0$ there exists a $0< T < \tau^x_0$ such that $X^x_t< \varepsilon$ for $t> T$. As $a= - \psi(0)= 0$ and $\psi'(0)\ge 0$, by Proposition 3.1 we see $\lim_{t\to\infty}X_t=0$. By Lemma~\ref{l5.1} we have $\lim_{y\uparrow x} \tau^y_0= \tau^x_0$. It follows that $\lim_{y\uparrow x}\lim_{n\rightarrow \infty}\zeta^{x,y}_n= \tau^x_0$. Then for sufficiently large $k,n$ and $y<x$ we have $k\land \zeta^{x,y}_n> T$, and thus
\beqnn
\sup_{t\ge 0}|X^x_t- X^y_t|
\ar\leq\ar
\sup_{0\leq s\leq k\land\zeta^{x,y}_n} |X^x_{s-}- X^y_{s-}|+\varepsilon\\
\ar=\ar
\sup_{0\leq s\leq k\land\zeta^{x,y}_n} |\xi_n^x(s-)- \xi_n^y(s-)|+\varepsilon.
\eeqnn
Since $\varepsilon>0$ can be arbitrarily small, we have a.s.\ $\lim_{y\uparrow x}\sup_{t\ge 0}|X^x_t- X^y_t|= 0$. \qed

\medskip

\noindent\emph{Proof of Corollary}~\ref{c3}.~ It follows from (\ref{1.5}), (\ref{1.4}) and dominated convergence that $\lim_{y\to\infty}$ $\mbf{E}(\tau^\infty_y)$ $= 0$. Then we have a.s.\ $\lim_{y\to\infty} \tau^\infty_y= 0$. \qed

In the sequel of this section, we assume $a= -\psi(0)= 0$, $\psi'(0)\geq 0$ and (\ref{1.4}) holds. By Theorem~\ref{t3.2} and Corollary~\ref{c3} we have a.s.\ $\tau^{\infty}_y\in (0,\infty)$ for $y\in (0,\infty)$ and $\lim_{y\to \infty} \tau^\infty_y= 0$. Now fix $y\in (0,\infty)$. For any positive random variable $\xi$ measurable with respect to $\mcr{F}_{\tau^\infty_y}$, consider the stochastic integral equation
\beqlb\label{e0.3}
X_t\ar=\ar \xi +\sqrt{2c}\int^{t}_0 X_s^{\theta/2}dB_{\tau^\infty_y+s} + \int_0^t\int_{(0,1]} \int^{X_{s-}^\theta}_0 z \tilde{M}(\tau^\infty_y+ds,dz,du)\cr
 \ar\ar\quad
-\, b \int^{t}_0 X_s^\theta ds + \int_0^{t}\int_{(1,\infty)}\int^{X_{s-}^\theta}_0 z M(\tau^\infty_y+ds,dz,du).
 \eeqlb

\blemma\label{p4}
As $k\to \infty$, in the supremum norm $(X^k_{\tau^\infty_y+t})_{t\ge 0}$ converges a.s.\ to a process $(X^{(y)}_t)_{t\ge 0}$, which is the pathwise unique solution to {\rm(\ref{e0.3})} with $\xi=y$.
\elemma

\proof Recall that $(X^k_t)_{t\ge 0}$ is the solution to (\ref{e0.1}) with $X^k_0= k$. It follows that $(Y_t)_{t\ge 0}:= (X^k_{\tau^\infty_y + t})_{t\ge 0}$ solves (\ref{e0.3}) with $Y_0= X^k_{\tau^\infty_y}$. By Theorem~\ref{t1.1}, there is a pathwise unique solution $(X^{(y)}_t)_{t\ge 0}$ to (\ref{e0.3}) with $\xi= y$. Using the strong Markov property we see $(X^k_{\tau^k_y+t})_{t\ge 0}= (X^y_t)_{t\ge 0}$ in distribution, where $(X^y_t)_{t\ge 0}$ is the solution to (\ref{e0.1}) with $X^y_0=y$. For any $\varepsilon>0$ we can choose $\delta>0$ so that $\mbf{P}(\sup_{0\le t\le \delta}|X^y_t-y|> \varepsilon)\le \varepsilon/2$. Since $\tau^k_y\uparrow \tau^\infty_y$ as $k\uparrow \infty$, there is $k_0\ge 1$ so that $\mbf{P}(\tau^\infty_y> \tau^k_y+\delta)\le \varepsilon/2$ for $k\ge k_0$. For $k\ge k_0$ we have
 \beqnn
\mbf{P}(|X^k_{\tau^\infty_y}-y|> \varepsilon)
 \ar\le\ar
\mbf{P}(|X^k_{\tau^\infty_y}-y|> \varepsilon, \tau^k_y\le \tau^\infty_y\le \tau^k_y+ \delta) + \varepsilon/2 \cr
 \ar\le\ar
\mbf{P}\Big(\sup_{0\le t\le \delta}|X^k_{\tau^k_y+t}-n|> \varepsilon\Big) + \varepsilon/2 \cr
 \ar\le\ar
\mbf{P}\Big(\sup_{0\le t\le \delta}|X^y_t-y|> \varepsilon\Big) + \varepsilon/2\le \varepsilon.
 \eeqnn
It follows that $\lim_{k\to\infty}X^k_{\tau^\infty_y}= y$ in probability. But $k\mapsto X^k_{\tau^\infty_y}$ is increasing by Theorem~\ref{t1.1}, so we also have a.s.\ $\lim_{k\to\infty}X^k_{\tau^\infty_y}= y$. Then the desired result follows from Proposition~\ref{l1.2}.\qed

For $y\in [0,\infty)$ let $(X^{(y)}_t)_{t\ge 0}$ be given by Lemma~\ref{p4}. By the pathwise uniqueness, for any $y\le x\in (0,\infty)$ and $t\ge \tau^\infty_y$ we have a.s.\ $X^{(x)}_{t-\tau^\infty_x} = X^{(y)}_{t-\tau^\infty_y}$. Then we can construct a positive c\`adl\`ag process $X^\infty= (X^\infty_t)_{t>0}$ such that $X^\infty_t= X^{(n)}_{t-\tau^\infty_n}$ a.s.\ for each $t\ge \tau^\infty_n$ and $n\ge 1$.

\bproposition\label{p5}
The process $X^\infty= (X^{\infty}_t)_{t> 0}$ defined above is a solution to {\rm(\ref{e0.1})} with initial state $\infty$. Moreover, we have a.s.\ $X^\infty_t\ge X^x_t$ for $t>0$ and $x\in (0,\infty)$ and $\tau^\infty_y = \inf\{t>0: X^\infty_t=y\}$ for $y\in (0,\infty)$.
\eproposition

\proof From the construction of $X^\infty$ it is clear that $X^\infty_{\tau^\infty_y}= X^{(y)}_0= y$ and $X^\infty$ satisfies (\ref{ne1}) for $t>r>0$. For any $t>0$ we can choose $n\ge 1$ so that $0< \tau^\infty_n\le t$. By Lemma~\ref{p4} we see $X^\infty_t= X^{(n)}_{t-\tau^\infty_n}\ge X^k_t$ for any $k\ge 1$. Then $X^\infty_t\ge X^x_t$ for any $t>0$ and $x\in (0,\infty)$. In particular, we get $X^\infty_t\ge X^x_{t}> y$ for $0< t< \tau^x_y$ if $y\le x\in (0,\infty)$. Since $\tau^x_y\uparrow \tau^\infty_y$ as $x\uparrow \infty$, we see $X^\infty_t> y$ for $0< t< \tau^\infty_y$, implying $\lim_{t\downarrow 0}X^\infty_t= \infty$ and $\tau^\infty_y= \inf\{t>0: X^\infty_t= y\}$. Then $X^\infty$ is a solution to {\rm(\ref{e0.1})} with initial state $\infty$.  \qed

\medskip

\noindent\emph{Proof of Theorem}~\ref{t3.3}.~ The existence of the solution $X^\infty= (X^{\infty}_t)_{t> 0}$ to (\ref{e0.1}) with initial state $\infty$ follows by Proposition~\ref{p5}. Fix $n\ge 1$. Since $\lim_{k\to\infty}\tau^k_{n}= \tau^\infty_{n}> \tau^\infty_{2n}$, we can choose sufficiently large $k\ge 1$ so that $\tau^k_{n}> \tau^\infty_{2n}$. Then for any $x\ge k$ we have $\tau^x_{n}> \tau^\infty_{2n}$, and so
\beqnn
\rho_\infty(X^x, X^\infty)
\le
\sup_{0\le t<\tau^x_n}(e^{-X^x_t}-e^{-X^\infty_t})\vee\sup_{t\ge \tau^\infty_{2n}} (X^\infty_t-X^x_t).
\eeqnn
From the construction of $X^\infty$ we see
 \beqnn
\sup_{t\ge\tau^\infty_{2n}}(X^\infty_t - X^x_t)
 =
\sup_{t\ge0}(X^\infty_{\tau^\infty_{2n}+t} - X^x_{\tau^\infty_{2n}+t})
 =
\sup_{t\ge0}(X^{(2n)}_{t} - X^x_{\tau^\infty_{2n}+t}).
 \eeqnn
By Lemma~\ref{p4}, the right-hand side vanishes as $x\to \infty$. Using Proposition~\ref{p5} we can see $X^\infty_t\ge X^x_t> n$ for $0< t< \tau^x_n$, and hence $\sup_{0< t< \tau^x_n} (e^{-X^x_t}-e^{-X^\infty_t})\le e^{-n}$. It follows that $\limsup_{x\to \infty} \rho_\infty(X^x, X^\infty)\le e^{-n}$. Since $n\ge 1$ can be arbitrary, we conclude $\lim_{x\to \infty} \rho_\infty(X^x, X^\infty)= 0$. Now suppose that $Y= (Y_t)_{t> 0}$ is another solution to (\ref{e0.1}) with initial state $\infty$. By right-continuity and comparison property we see a.s.\ $X^\infty_t= \lim_{x\to \infty} X^x_t\le Y_t$ for any $t>0$. For $n\ge 1$ let $(X_{n,t}^x)_{t> 1/n}$ be the pathwise unique solution to
\beqnn
X_{t}\ar=\ar x+\sqrt{2c}\int^{t}_{1/n} X_{s}^{\theta/2}dB_s + \int_{1/n}^t\int_{(0,1]} \int^{X_{s-}^\theta}_0 z \tilde{M}(ds,dz,du)\cr
 \ar\ar\quad
-\, b \int^{t}_{1/n} X_{s}^\theta ds + \int_{1/n}^{t}\int_{(1,\infty)} \int^{X_{s-}^\theta}_0 z M(ds,dz,du).
 \eeqnn
Let $X^\infty_{n,t}= \lim_{x\to\infty} X_{n,t}^x$ for $t\ge 1/n$. By the comparison property we see a.s.\ $Y_{t}\le X^\infty_{n,t}$ for $t\ge 1/n$. Notice that $X^\infty_{n,t}= X^\infty_{t-1/n}$ in distribution for each $t\ge 1/n$. Then $X^\infty_{n,t}\to X^\infty_t$ in distribution as $n\to\infty$. From the a.s.\ relation $X^\infty_t\le Y_{t}\le X^\infty_{n,t}$, it follows that $Y_t= X^\infty_t$ in distribution for $t>0$. Then we must have $Y_t=X^\infty_t$ a.s.\ for each $t>0$. That gives the pathwise uniqueness for (\ref{e0.1}) with initial state $\infty$. \qed

\section{Convergence of discrete-state processes}

\setcounter{equation}{0}

In this section, we study the convergence of rescaled discrete-state polynomial branching processes to continuous-state ones. Let us consider a sequence of generating functions $g_n, n=1,2,\ldots$ given by
 \beqlb\label{S2ax}
g_n(s) =\sum^\infty_{i=0} b_i^{(n)} s^i, \qquad s\in [0,1],
 \eeqlb
where $\{b_i^{(n)}: i= 0,1,\ldots,\infty\}$ is a discrete probability distribution. Let $\{\gamma_n: n=1,2,\ldots\}$ be a sequence of positive numbers. We define the two sequences of functions $\{\psi_n\}$ and $\{\phi_n\}$ by
 \beqlb\label{S1}
\phi_n(\lambda)=\gamma_n[g_n(e^{-\lambda/n})-e^{-\lambda/n}], \qquad \lambda\geq 0
 \eeqlb
and
 \beqlb\label{S2}
\psi_n(\lambda)=\gamma_n[g_n(1-\lambda/n)-(1-\lambda/n)],\qquad 0\leq\lambda\leq n.
 \eeqlb

\bproposition\label{tch01.3.5}
The sequence $\{\phi_n\}$ defined by {\rm(\ref{S1})} is Lipschitz uniformly on each bounded interval $[\alpha,\beta]\subset (0,\infty)$ if and only if so is the sequence $\{\psi_n\}$ defined by {\rm(\ref{S2})}. In this case, we have $\lim_{n\to\infty} |\psi_n(\lambda)-\phi_n(\lambda)| = 0$ uniformly on each bounded interval $[\alpha,\beta]\subset (0,\infty)$.
\eproposition

\proof Clearly, the sequences $\{\phi_n\}$ or $\{\psi_n\}$ are Lipschitz uniformly on some interval $[\alpha,\beta]\subset (0,\infty)$ if and only if the sequences of derivatives $\{\phi_n^\prime\}$ or $\{\psi_n^\prime\}$ are bounded uniformly on the interval. From (\ref{S1}) and (\ref{S2}) we have
 \beqnn
\phi_n^\prime(\lambda)
 =
n^{-1}\gamma_ne^{-\lambda/n}[1-g^\prime_n(e^{-\lambda/n})], \qquad \lambda\ge 0
 \eeqnn
and
 \beqnn
\psi^\prime_n(\lambda)
 =
n^{-1}\gamma_n[1-g^\prime_n(1-\lambda/n)], \qquad 0\le \lambda\le n.
 \eeqnn
Then $\{\phi_n^\prime\}$ is uniformly bounded on each bounded interval $[\alpha,\beta]\subset (0,\infty)$ if and only if so is $\{\psi_n^\prime\}$. That proves the first assertion. We next assume $\{\phi_n\}$ is Lipschitz uniformly on each bounded interval $[\alpha,\beta]\subset (0,\infty)$. Observe that
 \beqnn
\phi_n(\lambda) - \psi_n(\lambda)
 =
\gamma_n\big[g_n(e^{-\lambda/n}) - e^{-\lambda/n} - g_n(1-\lambda/n) + (1-\lambda/n)\big].
 \eeqnn
By the mean-value theorem, for $n\ge \beta$ and $\alpha\le \lambda\le \beta$ we have
 \beqlb\label{ch01.3.10}
\phi_n(\lambda)-\psi_n(\lambda)
 =
\gamma_n[g^\prime_n(\eta_n)-1](e^{-\lambda/n}-1+\lambda/n),
 \eeqlb
where $1-\lambda/n\le \eta_n := \eta_n(\lambda)\le e^{-\lambda/n}$. Choose sufficiently large $n_0\ge \beta$ so that $e^{-2\beta/n_0}\le 1-\beta/n_0$. For $n\ge n_0$ we have $e^{-2\beta/n}\le 1-\beta/n\le 1-\lambda/n$. It follows that $e^{-2\beta/n}\le \eta_n\le e^{-\alpha/n}$ for $\alpha\le \lambda\le \beta$. By the monotonicity of $z\mapsto g^\prime(z)$,
 \beqnn
n^{-1}\gamma_n |g^\prime_n(\eta_n)-1|
 \le
\sup_{\alpha\le \lambda\le 2\beta} n^{-1}\gamma_n|g^\prime_n(e^{-\lambda/n})-1|
 =
\sup_{\alpha\le \lambda\le 2\beta} e^{\lambda/n}|\phi_n^\prime(\lambda)|.
 \eeqnn
Then $\{n^{-1}\gamma_n|g^\prime_n(\eta_n)-1|: n\ge n_0\}$ is a bounded sequence. Since $\lim_{n\to \infty}n(e^{-\lambda/n}-1+\lambda/n) = 0$ uniformly on $[\alpha,\beta]$, the desired result follows by (\ref{ch01.3.10}). \qed

\bproposition\label{p3.6}
For any function $\psi$ on $[0,\infty)$ with representation {\rm(\ref{0.1})} there is
a sequence $\{\phi_n\}$ in form {\rm(\ref{S1})} so that $\lim_{n\to\infty} \phi_n(\lambda) = \psi(\lambda)$ for $\lambda\geq 0$.
\eproposition

\proof By Proposition~\ref{tch01.3.5} it is sufficient to construct a sequence $\{\psi_n\}$ in form (\ref{S2}) that is Lipschitz uniformly on $[\alpha,\beta]$ and $\lim_{n\to \infty} \psi_n(\lambda) = \psi(\lambda)$ uniformly on $[0,\beta]$ for any $\beta>\alpha>0$. In view of (\ref{0.1}), we can write
 \beqnn
\psi(\lambda) = -a + b_n\lambda + c\lambda^2 + \int_{(0,\infty)} \big(e^{-\lambda
u} - 1 + \lambda u1_{\{u\le \sqrt{n}\}}\big) m(du),
 \eeqnn
where
 \beqnn
b_n = b - \int_{(1,\sqrt{n}]} u m(du).
 \eeqnn
Observe that $|b_n|\le |b| + m(1,\infty)\sqrt{n}$. Let $\gamma_{1,n} = n$ and $g_{1,n}(z) = (1-n^{-2}a)z$. Let $\psi_{1,n}(\lambda)$ be defined by (\ref{S2}) with $(\gamma_n,g_n)$ replaced by $(\gamma_{1,n},g_{1,n})$. Then we have $\psi_{1,n}(\lambda) = -a(1-\lambda/n)$. Following the proof of Proposition~4.4 in Li (2011, p.93) one can find a sequence of positive numbers $\{\alpha_{2,n}\}$ and a sequence of probability generating functions $\{g_{2,n}\}$ so that the function $\psi_{2,n}(\lambda)$ defined by (\ref{S2}) from $(\alpha_{2,n}, g_{2,n})$ is given by
 \beqnn
\psi_{2,n}(\lambda) = b_n\lambda + \frac{1}{2n}(|b_n|-b_n)\lambda^2 + c\lambda^2 + \int_{(0,\sqrt{n}]} \big(e^{-\lambda u} - 1 + \lambda u\big) m(du).
 \eeqnn
Let $\gamma_n = \gamma_{1,n} + \gamma_{2,n}$ and $g_n(z) = \gamma_n^{-1} [\gamma_{1,n}g_{1,n}(z) + \gamma_{2,n}g_{2,n}(z)]$. Then the sequence $\{\psi_n(\lambda)\}$ defined by (\ref{S2}) is equal to $\{\psi_{1,n}(\lambda) + \psi_{2,n}(\lambda)\}$, which clearly possesses the required properties. \qed

\medskip

\noindent\emph{Proof of Theorem}~\ref{t3.1}.~ By Proposition~\ref{t1.0'}, the generalized Lamperti transform $Y=L_\theta(X)$ is a L\'{e}vy process with Laplace exponent $-\psi$ stopped at $0$. By Proposition \ref{p3.6} there is a sequence $\{\phi_n\}$ in form (\ref{S1}) so that $\lim_{n\to\infty} \phi_n(\lambda) = \psi(\lambda)$ for $\lambda\geq 0$. By adjusting the parameters, we may assume the probability distribution $\{b_i^{(n)}: i=0,1,2,\ldots,\infty\}$ satisfies $b_1^{(n)}=0$. Let $Z_n= (Z_n(t): t\ge 0)$ be a compound Poisson process on the state space $\{0,\pm 1,\pm 2,\ldots,\infty\}$ with $Q$-matrix defined by
 \beqnn
\rho_n(i,j) = \left\{
 \begin{array}{lcl}
 b^{(n)}_{j-i+1}, && {i+1\le j< \infty,} \cr
 -1, && {i= j< \infty,} \cr
 b^{(n)}_0, && {i-1= j< \infty,} \cr
 0, && \mbox{otherwise.}
 \end{array} \right.
 \eeqnn
Then $Z_n$ has Laplace exponent $e^{\lambda}[e^{-\lambda} - g_n(e^{-\lambda})]$. Let $T_n = \inf\{t\ge 0: Z_n(t) = 0\}$ and let $Y_n = (Z_n(t\land T_n): t\ge 0)$ be the stopped process. Set $Z^{(n)}(t) = n^{-1}Z_n(\gamma_nt)$. The rescaled compound Poisson process $Z^{(n)} = (Z^{(n)}(t): t\ge 0)$ has Laplace exponent
 $$
e^{\lambda/n}\phi_n(\lambda) = n\gamma_n e^{\lambda/n}[e^{-\lambda/n}-g_n(e^{-\lambda/n})].
 $$
Let $T^{(n)} = \inf\{t\ge 0: Z^{(n)}(t) = 0\}$ and let $Y^{(n)} = (Z^{(n)}(t\land T^{(n)}): t\ge 0)$ be the stopped process. The inverse generalized Lamperti transforms $X_n:= J_\theta(Y_n)$ and $X^{(n)}:= J_\theta(Y^{(n)})$ can be defined similarly as in the introduction. By a simple extension of Theorem 2.1 in Chen et al.\ (2008), one can see $X_n$ is a discrete-state polynomial branching process with $Q$-matrix given by
 \beqnn
q_n(i,j) = \left\{
 \begin{array}{lcl}
 i^\theta b^{(n)}_{j-i+1}, && {i+1\leq j< \infty, i\ge 1,} \cr
 -i^\theta, && {1\le i= j< \infty,}\cr
 i^\theta b^{(n)}_0, && {0\le i-1= j< \infty,} \cr
 0, && \mbox{otherwise.}
 \end{array} \right.
 \eeqnn
Then $X^{(n)}$ is a rescaled discrete-state polynomial branching process. Since
 $$
\lim_{n\to\infty} e^{\lambda/n}\phi_n(\lambda)
 =
\lim_{n\to\infty} \phi_n(\lambda) = \psi(\lambda),
 \qquad \lambda\geq 0,
 $$
by Proposition 6 in Caballero et al.\ (2009) we see ${Y^{(n)}}\to Y$ weakly in $(D,d_\infty)$. By a slight generalization of Proposition 5 in Caballero et al.\ (2009), one can see the transformation $J_\theta$ is continuous on $(D,d_{\infty})$. Then $X^{(n)}\to X = J_\theta(Y)$ in $(D,d_\infty)$. \qed

\bigskip\bigskip

\noindent\textbf{Acknowledgments.}~ I am very grateful to Professors Mu-Fa Chen and Yong-Hua Mao for their supervision and encouragement. I thank the associate editor and the two anonymous referees for their thoughtful comments which have helped me improving the presentation of the results. I also thank Professors Zenghu Li and Xiaowen Zhou for their helpful comments. This work is supported by NSERC (RGPIN-2016-06704) and NSFC (No.\,11771046).

\end{document}